\documentclass[11pt,a4paper]{article}

\usepackage[T1]{fontenc}
\usepackage[utf8]{inputenc}
\usepackage{lmodern}
\usepackage{amsmath,amssymb,amsthm,mathtools}
\IfFileExists{stmaryrd.sty}{\usepackage{stmaryrd}}{%
  \newcommand{\llbracket}{[\![}%
  \newcommand{\rrbracket}{]\!]}%
}
\usepackage[a4paper,margin=29mm]{geometry}
\usepackage{microtype}

\makeatletter
\renewcommand\subsection{%
  \@startsection{subsection}{2}{\z@}%
    {1.6ex \@plus .6ex \@minus .2ex}%
    {-1em}%
    {\normalfont\normalsize\bfseries}%
}
\makeatother
\usepackage[hidelinks]{hyperref}

\hypersetup{
  pdftitle={Euclidean isoperimetric inequalities in CAT(0)},
  pdfsubject={Notes on sharp integral fillings and their current-theoretic preliminaries}
}

\newtheorem{theorem}{Theorem}[section]
\newtheorem{proposition}[theorem]{Proposition}
\newtheorem{lemma}[theorem]{Lemma}
\newtheorem{corollary}[theorem]{Corollary}
\numberwithin{equation}{section}
\newtheorem{remark}[theorem]{Remark}

\newcommand{\R}{\mathbb R}

\newcommand{\I}{\mathbf I}
\newcommand{\M}{\mathbf M}
\newcommand{\HH}{\mathcal H}
\newcommand{\LL}{\mathcal L}
\newcommand{\CAT}{\operatorname{CAT}}
\DeclareMathOperator{\spt}{spt}
\DeclareMathOperator{\md}{md}
\DeclareMathOperator{\Lip}{Lip}

\title{Euclidean filling inequalities in CAT(0) spaces}
\author{Alexander Lytchak}
\date{}

\begin{document}
\maketitle
\vspace{-1.4em}

\begin{abstract}
We clarify and streamline OpenAI's proof of the sharp Euclidean filling inequalities in CAT(0) spaces. We extend the result to complete spaces without assuming properness or compact support,  and we establish rigidity in the equality case.
\end{abstract}

\section{Introduction}

\subsection{Results}
We establish the optimal Euclidean filling inequality for integral cycles
in complete $\CAT(0)$ spaces, without assuming properness of the ambient
space or compactness of the support. Let $\omega_n$ denote the volume of
the Euclidean unit ball in $\R^n$, let $\I_n(Y)$ be the group of integral
$n$-currents in $Y$, and let $\M(A)$ denote the mass of an integral
current $A$. For an integral $n$-cycle $A\in\I_n(Y)$, its filling volume
in $Y$ is
\[
 V(A)=V_Y(A)=\inf\{\M(B):B\in\I_{n+1}(Y),\ \partial B=A\}.
\]

\begin{theorem}\label{thm:main}
Let $Y$ be a complete $\CAT(0)$ space and let $n$ be a positive integer.
Every integral current $A\in\I_n(Y)$ with $\partial A=0$ admits an
integral filling $B\in\I_{n+1}(Y)$ satisfying
\[
  \partial B=A,
  \qquad
  (n+1)^{n+1}\omega_{n+1}\,\M(B)^n\leq\M(A)^{n+1}.
\]
\end{theorem}

In terms of filling volume, the inequality reads
\begin{equation} \label{eq:main}
  V(A)\leq c_n\M(A)^{(n+1)/n},
  \qquad
  c_n=\frac{1}{(n+1)^{(n+1)/n}\omega_{n+1}^{1/n}}.
\end{equation}

The  equality case is
characterized by the following rigidity statement.

\begin{theorem}\label{thm:rigidity}
Let $Y$ be a complete $\CAT(0)$ space and let $n\geq1$. For a nonzero $n$-cycle $A\in\I_n(Y)$, equality holds in \eqref{eq:main} if
and only if there are $R>0$ and an isometric embedding
$
 \iota:\overline B^{\,n+1}(0,R)\longrightarrow Y$
of the Euclidean ball such that
\[
 A=\iota_\#\partial\llbracket \overline B^{n+1}(0,R)\rrbracket.
\]
\end{theorem}

\subsection{History}
The proof of Theorem~\ref{thm:main} is a modification and simplification
of the proof by OpenAI \cite{OAI26}, where the inequality is established
for compactly supported cycles in proper $\CAT(0)$ spaces in dimensions
$n\geq2$. The main predecessors are 
\begin{itemize}
\item Reshetnyak's  majorization theorem  \cite{Reshetnyak}, implying Theorem \ref{thm:main} for $n=1$;
\item Almgren's proof of both theorems above in the case $Y=\R^N$ \cite{Almgren};
\item  Proof of  \eqref{eq:main} with a non-optimal constant $D_n>c_n$ by  Gromov  \cite{Gromov} on Cartan--Hadamard manifolds;
\item Proof of  \eqref{eq:main} with a non-optimal constant $D_n>c_n$ in  the general  CAT(0) setting  by  Wenger \cite{W05};
\item  Proof of the existence of a minimal filling of any $n$-cycle by Ambrosio--Kirchheim and Wenger \cite{AK00}, \cite{W14};
\item Proof   of the optimal filling inequality and the rigidity case in Cartan--Hadamard  manifolds  with $n=2$ by Schulze \cite{Schulze}. 
\end{itemize}

Besides, there is a rich history of results for  smooth ambient manifold $Y$ of dimension $n+1$, in which case the result has been known as the Cartan--Hadamard conjecture.   We mention Weil's result \cite{Weil}  in dimension $n+1=2$, Croke's \cite{Croke} in dimension $n+1=4$, and Kleiner's \cite{Kleiner}  in dimension $n+1=3$. The recent AI era in mathematics has seen an explosion of results in this area: Chen--Ghomi--Wang \cite{Ghomi3,Ghomi2} resolved the cases $n+1\leq9$ (September 10 and 17, 2026). OpenAI's manuscript \cite{OAI26+} resolving the general case  was publicly released on October 6, 2026. Agnoletto--Da Silva--Nardulli--Resende published  a different proof of the Cartan--Hadamard conjecture on October 27, 2026 \cite{Agnoletto}, building on some ideas from \cite{OAI26+}.

The recent paper  \cite{Agnoletto} proves the natural generalization of the Cartan--Hadamard conjecture also for upper curvature bounds $\kappa <0$ together with corresponding rigidity statements.  It is  a natural question whether  the main results of the present paper are generalizable  to CAT($\kappa$)-spaces for $\kappa\neq 0$.

\subsection{Main steps}
In Section \ref{sec:prelim} we recall preliminaries on integral currents.  In Section \ref{sec:ekeland}, we prove, following the classical approach in \cite{Almgren}, that Theorem \ref{thm:main}  is equivalent to the statement (in each dimension $n$ separately):
\begin{proposition}\label{prop:variational}
Let $Y$ be a complete $\CAT(0)$ space, let $n$ be a positive integer, and
let $A\in\I_n(Y)$ be a nonzero cycle. Suppose that, for some $c>0$,
\begin{equation} \label{eq:power}
  c\cdot \bigl(\M(A)^{(n+1)/n}-\M(B)^{(n+1)/n}\bigr)
  \leq V(A-B)
\end{equation}
for every integral $n$-cycle $B\in\I_n(Y)$. Then $c\leq c_n$.
\end{proposition}
If $Y$ is compact, one can take as $A$  a maximizer of the  functional  $B\to V(B)-c \M (B)^{\frac {n+1} n} $.  If $Y$ is non-compact, Ekeland's variational principle \cite{Ek74} provides the required cycle $A$.

In Section  \ref{sec:compact-support} we start the proof of Theorem \ref{thm:main} and Proposition \ref{prop:variational} by induction on $n$.  The case $n=1$ is easily deduced from Reshetnyak's majorization theorem.  For the induction step, assume that the results are true in dimensions at most $n-1$. Assuming that Theorem \ref{thm:main} fails, we fix an \emph{extremal cycle} $A$
as in Proposition \ref{prop:variational}, with a super-Euclidean constant $c>c_n$.  In many  subsequent steps we investigate the structure of this cycle and prove finally that it cannot exist.

First, we choose an appropriate rescaling of $Y$, so that $m=\M(A)$ is smaller than the volume $s_n$ of the Euclidean $n$-sphere and such that  the linearized version of \eqref{eq:power} becomes
\begin{equation} \label{eq:linear}
	\M(A)-\M(B)
	\leq n\cdot V(A-B)\,.
\end{equation}
whenever $V(A-B)<\delta$, for  sufficiently small  $\delta >0$.  

This has strong regularity implications on $A$, its support $K$ and its mass measure $\mu=\|A\|$.  First, \eqref{eq:linear}, the optimal filling inequality in dimension $n-1$ 
and   the (non-optimal) Euclidean-type filling  inequality  in dimension $n$  provided by \cite{W05} imply that $(K,\|A\|)$ satisfies (in a precise sense) an almost Euclidean isoperimetric inequality  for small masses.   
Then  a standard argument shows that small balls in $K$ have a uniform lower bound on their  $\mu$-volume.   This leads to compactness of $K$ and to  lower $n$-density  $\theta_y \geq 1$  at every $y\in K$.

In Section \ref{sec:radial} we estimate mass changes of integral currents  along  radial contractions towards an arbitrary center $y\in Y$.  These elementary CAT(0) estimates together with \eqref{eq:linear} amount to the following conclusion, valid for every nonnegative Lipschitz function $g:K\to \R$ and for the distance function $r(x)=d(y,x)$:
\begin{equation}  \label{eq:mainineq}
	\int\bigl(ng+r\langle\nabla^A r,\nabla^A g\rangle\bigr)\,d\mu
	\leq n\int rg\sqrt{1-|\nabla^A r|^2}\,d\mu,
\end{equation}
where  $\nabla ^A g$ denotes the $\mu$-almost everywhere defined gradient  of  the Lipschitz function $g$ on the rectifiable set $K$. 

The inequality \eqref{eq:mainineq} has a natural geometric
interpretation. For closed smooth submanifolds of Hadamard manifolds,
its validity for all centers and nonnegative test functions
is equivalent to a uniform bound $|H|\leq n$ on the mean-curvature vector.
A similar interpretation is possible for integral currents
in smooth ambient manifolds, using generalized mean curvature
in the sense of varifolds. Thus, the subsequent argument also
provides sharp Sobolev and volume estimates for submanifolds
with bounded mean curvature. Such estimates play an important
role in the proofs of the Cartan--Hadamard conjecture.
Related, more general Sobolev-type inequalities for
submanifolds of arbitrary codimension have recently been
obtained by Chen and Ma \cite{ChenMa}.

The remaining task is to find an appropriate center $y\in Y$ and a non-negative function $g$, contradicting \eqref{eq:mainineq}. In Section \ref{sec:dimension-two} we verify that for  $n=2$, we can  choose any center $y\in K$ and consider  test functions \[g_{\varepsilon}=\frac 1 {\max \{r^{2}, \varepsilon ^{2}\}}\,\;\; \; \; \varepsilon >0,    \]   
used in a related  setting   in \cite{Schulze}.  For $\varepsilon \downarrow 0$ the inequality \eqref{eq:mainineq} will imply $m\geq 4\pi=s_2$  contradicting  the previously fixed  normalization.

As simple as the choice of  appropriate centers and  test functions $g$ for \eqref{eq:mainineq} is in dimension two,  so complicated  does it appear to be   in higher dimensions. 
The solution proposed in \cite{OAI26} for all $n\geq 3$  appears very mysterious to me. It uses  a family of test functions derived from a solution $v$ of the Yamabe-type problem on the integral current $A$.
  Unfortunately, I could not find any a priori reason why this might lead to a solution.  A mathematician with a better analytic background may hopefully identify some reasons for choosing this path.  I just describe the details below, being unable to explain it from a higher perspective.

Sections \ref{sec:sobolev-spaces}, \ref{sec:intrinsic-sobolev}, \ref{sec:sobolev-minimizer} establish the existence of a function on $K$ solving a natural generalization of the Yamabe problem.  More precisely, Section \ref{sec:sobolev-spaces} recalls the basic theory of Sobolev functions on the metric measure space $(K,d, \mu)$ and verifies the Rellich compactness of the embedding of the Sobolev space $\mathcal H=W^{1,2}(K)$ into $L^q(\mu)$, for every $1\leq q<2$.  This part of the argument works  for any compactly supported integral current, not only for our extremal cycle  $A$.    Finally, we observe in Section \ref{sec:sobolev-spaces} that \eqref{eq:mainineq} extends to all nonnegative functions $g$ in the Sobolev space $\mathcal H=W^{1,2}(K,d,\mu)$.

In Section~\ref{sec:intrinsic-sobolev} we deduce from the
almost Euclidean isoperimetric inequality for small subsets
of $K$ an almost sharp Sobolev inequality with the critical
exponent $p=\frac{2n}{n-2}$ and an almost optimal Euclidean leading constant. This is an analogue of Aubin's almost sharp Sobolev
inequality for closed Riemannian manifolds~\cite{Aubin76}.
For such manifolds, the inequality with the exact
Euclidean leading constant was subsequently proved by
Hebey--Vaugon~\cite{HV96}.
 The results follow readily from the  P\'olya--Szeg\H{o} rearrangement theorem
 of Nobili--Violo~\cite{NV25} and the sharp Euclidean Sobolev inequality
of Aubin and Talenti~\cite{Aubin76,Talenti}.

In Section~\ref{sec:sobolev-minimizer}, we apply  Aubin's approach to the Yamabe problem in the generality developed by Akutagawa--Carron--Mazzeo~\cite{ACM19}. Using the 
almost sharp Sobolev inequality and the mass
assumption $m<s_n$, this approach provides a  bounded nonnegative minimizer in $\mathcal H$  of a natural Sobolev quotient.  After rescaling, this  minimizer  weakly solves
\[
-\frac4{n-2}\Delta_A v+nv
=nv^{\frac{n+2}{n-2}},
\]
  where $\Delta_A$ denotes the  Laplacian on $(K,d,\mu)$ associated with the Cheeger energy. This equation is the metric-measure analog of the Yamabe equation on the unit sphere.

In Section \ref{sec10-ce-section}, we apply the variational inequality
\eqref{eq:mainineq} with $g=v^p$ and use a CAT(0) variance inequality
to find a center $y\in K$ with a sufficiently large weighted second
moment.

In Section~\ref{sec:energy-monotonicity},  we
construct a one-parameter family $v_a$, $a\geq0$, starting at the
Sobolev minimizer $v_0=v$. While  \(v=v_0\) is chosen independently of the center \(y\), the variation \(v_a\) depends on \(y\). 
  The construction of the variation bears a strong formal similarity with the Aubin--Talenti bubbles \cite{Aubin76}, \cite{Talenti}:  these  are functions on $\R^n$,  for which the optimal Sobolev inequality  is an equality.  A lengthy computation identifies a monotonic quantity along this variation, which leads to an upper bound for the energy $Q(v_a)$. This upper bound crucially  depends on the choice of our center $y$ in Section \ref{sec10-ce-section}.  Indeed, this choice controls the initial velocity of $Q(v_a)$, whose upper bound ultimately propagates to infinity.

For $a\to \infty$ the function $v_a$ concentrates  at  $y$ and approaches there an Aubin--Talenti bubble after rescaling. 
In Section~\ref{sec:concentration-completion}, we use this convergence to derive  a lower bound on the energies $Q(v_a)$ for large $a$,  essentially from the corresponding explicitly known values of energies of the Aubin--Talenti bubbles.
 For $a\to \infty$ this leads to a contradiction with the previous upper bounds on $Q(v_a)$.

Section~\ref{sec:equality-rigidity} is devoted to
 the equality case stated in Theorem \ref{thm:rigidity}. We follow    the steps explained above  and  analyze the equality case in each step.  A central role is played by  a result in CAT(0) geometry, which seems interesting enough to be stated in the introduction. It turns out that the equality case of the variance inequality for measures in CAT(0) spaces $Y$ \cite[Theorem 4.9(iii)]{Sturm} characterizes subsets of $Y$ whose convex hull is flat.
 
 \begin{proposition}\label{rig:variance-flatness}
 Let $\sigma$ be a Borel probability measure on a complete
 	 $\CAT(0)$ space $Y$.  Let the support $K$  of $\sigma$ be separable and assume  that   $I(y)=\int _Yd(x,y)^2\,d\sigma (x) $ is finite  for all $y\in Y$. Let $o$ denote  the 
 	barycenter of $\sigma$. If
 	\begin{equation}\label{rig:variance-flatness-assumption}
 		\int\int d(x,y)^2\,d\sigma(x)\,d\sigma(y)
 		=2\int d(o,x)^2\,d\sigma(x),
 	\end{equation}
 	then the closed convex hull of $K$ is isometric to a closed convex
 	subset of a Hilbert space.
 \end{proposition}
 On the other hand, for any  flat convex subset $K$ of $Y$ and any probability measure  $\sigma$ with finite second moment, which is supported on $K$, equality holds in \eqref{rig:variance-flatness-assumption}.

 \subsection{Acknowledgements}
 This paper arose in close collaboration with ChatGPT, an AI assistant
 developed by OpenAI. The collaboration involved exploring proof strategies,
 refining arguments, locating references, and drafting and revising the
 manuscript. The author takes full responsibility for the mathematical
 content and the final text.
 
 The author was supported in part by the German Research Foundation under Germany’s Excellence Strategy – EXC-2047/1 – 390685813.

\section{Notation and preliminaries} \label{sec:prelim}
The  Hausdorff measure $\HH^k$
is normalized to agree with Lebesgue measure $\mathcal L^k$ on $\R^k$.

Throughout, $Y$ is a complete $\CAT(0)$ space, not necessarily proper.
We write $B(x,r)$ for an open metric ball, $T_xY$ for the Alexandrov
tangent cone at $x\in Y$, and $\log_x:Y\to T_xY$ for the canonical
$1$-Lipschitz logarithm map.  Denote by  $d_x:Y\to \R$  the distance function to the point  $x\in Y$.

  For  $C \in \I_k(Y)$ we denote  by $\|C\|$ its mass measure, a finite Borel measure on $Y$.

\begin{proposition}\label{prop:wenger-fillings}
For each $k\geq1$ there is a constant $D_k<\infty$, depending only on
$k$, such that every integral $k$-cycle $T$ in $Y$ satisfies
\[
 V_Y(T)\leq D_k\M(T)^{(k+1)/k}.
\]
Moreover, the infimum defining $V_Y(T)$ is attained by an integral
$(k+1)$-current.
\end{proposition}

\begin{proof}
The inequality is \cite[Corollary~1.4]{W05}; the last statement  is \cite[Theorem~1.3]{W14}.
\end{proof}

The following result can be found directly after \cite[Corollary 4.4]{W06}.

\begin{lemma}\label{lem:lower-density}
For every $C\in\I_k(Y)$, $k\geq1$, the density is well-defined and at least $1$  on $Y$:
\[
  \lim_{r\downarrow0}\frac{\|C\|(B(y,r))}{\omega _kr^k}\geq 1
  \qquad\text{for $\|C\|$-almost every $y\in Y$.}
\]
\end{lemma}

We recall the following  version of Rademacher's theorem:

\begin{proposition} \label{prop:differential}
	Let $E\subset \R^k$ be a Borel subset, and let $\phi:E\to Y$ be a Lipschitz map. Then for almost every $z\in E$ there exists a well-defined differential $d\phi _z:T_z\R^k=\R^k \to T_x Y$, where $x=\phi (z)$ such that the following holds true:
	\begin{equation}\label{eq:chart-differentiability}
		d_{T_xY}\bigl(\log_x\phi(z+h),d\phi_z(h)\bigr)
		=o(|h|)\qquad(h\to0,\ z+h\in E);
	\end{equation}
	\begin{equation} \label{eq:Radem}
		d\phi _z( \R^k) \subset T_xY \; \; \text{is isometric to} \;\; \R^l \; \; \text{with} \; \; l\leq k  \; \;  \text{and}  \; \; d\phi _z:\R^k \to \R^l \; \; \text{linear}.
	\end{equation}
\end{proposition}

\begin{proof}
The existence and uniqueness of $d\phi _z$ with \eqref{eq:chart-differentiability} and \eqref{eq:Radem} is \cite[Theorem 1.6]{Lytchak}, see also \cite[Theorem 13.12]{AKP}.
\end{proof}

\begin{theorem}\label{thm:charts}
Let $k\geq1$ and $C\in\I_k(Y)$. For $\|C\|$-almost every $x\in Y$
there is a  Euclidean $k$-dimensional approximate tangent
space
\[
  T_xC\subset T_xY.
\]
For any $\varepsilon>0$ there are countably many compact sets
$E_i\subset\R^k$, $(1+\varepsilon)$-biLipschitz maps
$\phi_i:E_i\to Y$, with pairwise disjoint images
$Z_i=\phi_i(E_i)\subset\spt C$ and $\theta _i \in \mathbb N$  such that
\begin{equation}\label{eq:representation}
	C(b,\pi_1,\ldots,\pi_k)
	=\sum_i \theta _i \int_{E_i}b(\phi_i(z))
	\det D\bigl((\pi_1,\ldots,\pi_k)\circ\phi_i\bigr)(z)\,dz,
\end{equation}
for all Lipschitz functions $\pi_a$ and all bounded Lipschitz functions $b$.

At almost all  $z\in E_i$ the differential $d(\phi_i )_z$ defined as in Proposition \ref{prop:differential}
is a linear isomorphism from $\R^k$ to $T_xC$, where $x=\phi_i(z)$.

Denote by $J_i(z)$  the $k$-Jacobian  of $\phi_i$ at $z$, that is  the norm of the determinant of the linear isomorphism $d(\phi_i ) _z$ between $\R^k$ and the image $T_xC$.  Then
$\|C\|(Y\setminus\bigcup_i Z_i)=0$, and
\begin{equation}\label{eq:mass-formula}
	\|C\|=\sum_i  \theta _i  \cdot (\phi_i)_\#
	\bigl(J_i\,\LL^k\llcorner E_i\bigr).
\end{equation}
\end{theorem}

\begin{remark} Note that the metric differential $\md \phi_i (z)$ of $\phi _i$ at a point $z$ is  $\md\phi_i(z)(v)=|d\phi_i(z)(v)|$ and thus defined by a Euclidean norm.  The bilipschitz property implies that the Jacobians $J_i$ satisfy $|J_i(z)-1|\leq 2k\varepsilon$, at almost every $z$, once $\varepsilon$ is  small enough.
	\end{remark}

\begin{proof}
The  representation \eqref{eq:representation} with $\phi_i$ being bilipschitz but possibly not $(1+\varepsilon)$-bilipschitz  and $\theta _i$ possibly non-constant integrable, integer-valued 
functions  is 
\cite[Theorem~4.5]{AK00}.
 The metric differentials of maps from compact subsets  of $\R^k$  to
$\CAT(0)$ spaces are induced by inner products  as we have recalled above. 

Therefore, \cite[Lemma~2.12]{BT24} allows us to write each $Z_i$ as a disjoint union  of $(1+\varepsilon)$-bilipschitz images of  compact subsets of $\R^k$ and an $\HH^k$-negligible  set.  The proof of  \cite[Theorem~4.5]{AK00}  provides   \eqref{eq:representation} with $\phi _i$ now being $(1+\varepsilon)$-bilipschitz.

By further subdividing the chart domains $E_i$,  we may assume that the multiplicities are nonzero constants. Reversing the orientation of the parametrizations where necessary, we may further assume that all multiplicities are positive.

The statement $\|C\|(Y\setminus\bigcup_i Z_i)=0$ follows from 
\eqref{eq:representation}.

Fix now such $\phi _i$ and define $T_xC$ to be $d(\phi _i)_z (\R^k)$, whenever $x=\phi _i (z)\in Z_i$ and the differential of $\phi _i$ at the Lebesgue density point $z\in E_i$  is defined by Proposition \ref{prop:differential}.
Since $\phi_i$ is $(1+\varepsilon)$-bilipschitz the differential $d(\phi_i)_z$ is $(1+\varepsilon)$-bilipschitz as well, hence $T_xC$ is a $k$-dimensional Euclidean space.
The definition is independent of the charts: on the intersection of any two chart images, the transition map is bilipschitz and its differential is a linear isomorphism almost everywhere. The chain rule for the differentials into \(T_xY\) therefore identifies their images. Thus \(T_xC\) is well-defined for \(\|C\|\)-almost every \(x\).

Finally,   \eqref{eq:mass-formula} is
\cite[Theorem~9.5]{AK00} together with
\cite[Corollary~8]{K94}.
\end{proof}

\begin{proposition}\label{prop:tangential-differentials}
Let $C\in\I_k(Y)$, $k\geq1$, and let $Z$ be a complete
$\CAT(0)$ space. If $F:\spt C\to Z$ is Lipschitz, then, for
$\|C\|$-almost every $x$, there is a well-defined differential
\[
  d^C F_x:T_xC\longrightarrow T_{F(x)}Z
\]
which is a linear map onto a Euclidean subspace of $T_{F(x)}Z$.  
 For arbitrary  charts $\phi_i$   of Theorem~\ref{thm:charts} such that 
 $\phi _i (z)=x$  and $\phi _i$ is differentiable at $z$ as in Proposition \ref{prop:differential},  we have
\begin{equation}\label{eq:tangential-chain-rule}
  d^C F_x\circ d\phi_i(z)=d(F\circ\phi_i)_z,
  \qquad x=\phi_i(z),
\end{equation}
\end{proposition}

\begin{proof}
Fix a decomposition of $\spt C$ into  the sets $Z_i=\phi _i (E_i)$  and a negligible remainder as in Theorem \ref{thm:charts}.   Consider the set of Lebesgue density points $z\in E_i$, at which the differentials $d(\phi_i)_z$ and $d(F\circ \phi _i)_z$  as in Proposition \ref{prop:differential} exist.   Then, for $x=\phi_i (z)$, the differential   $d(\phi _i)_z:\R^k \to T_xC$  is $(1+\varepsilon)$-bilipschitz and $d^CF_x$  is well-defined through \eqref{eq:tangential-chain-rule}.

For this fixed representation of $\spt C$ as a union of charts, this  defines $d^CF_x:T_xC \to T_{F(x)}Z$ for $\|C\|$-almost every $x\in \spt C$.  
If another decomposition is given, the chain rule  shows independence of this differential at almost all points.
\end{proof}

 Let $F$ in Proposition \ref{prop:tangential-differentials} be $L$-Lipschitz in a neighborhood of $x$. If $\phi_i$ appearing in \eqref{eq:tangential-chain-rule} are $(1+\varepsilon)$-bilipschitz then $d^CF_x$ is $(1+2\varepsilon)^2\cdot L$-Lipschitz, due to \eqref{eq:tangential-chain-rule}.  Letting $\varepsilon$ go to $0$ we deduce that $d^CF_x$  is $L$-Lipschitz at $\|C\|$-almost every point $x$.

 Applying Proposition \ref{prop:tangential-differentials}  to Lipschitz functions $f:\spt C\to \R$  we obtain $\|C\|$-almost everywhere defined linear functionals $d^Cf_x:T_xC\to \R$, hence a  vector
 $\nabla^Cf (x) \in T_xC$ defined almost everywhere  with $|\nabla^C f(x)|\leq\Lip(f)$, and such that for all $v\in T_xC$:   
 \[
  d^C f_x(v)=\langle\nabla^C f(x),v\rangle.
  \]

\section{Application of Ekeland's principle}\label{sec:ekeland}
The following result is folklore, but we have been unable to find a precise reference.

\begin{lemma}\label{lem:filling-metric}
Let $n$ be a positive integer and let $Y$ be a complete $\CAT(0)$
space. On the group of integral $n$-cycles in $Y$ the formula
\[
 d_V(A,B):=V_Y(A-B)
\]
defines a finite-valued metric. For every $t>0$, the subspace
\[
 \mathcal Z_t:=\{A\in\I_n(Y):\partial A=0,\ \M(A)\le t\}
\]
is complete in this metric. The mass functional is lower
semicontinuous with respect to $d_V$.
\end{lemma}

\begin{proof}
Proposition~\ref{prop:wenger-fillings} makes $d_V$ finite. Symmetry
and the triangle inequality are immediate. If $V_Y(C)=0$, fillings
$S_j$ of $C$ with $\M(S_j)\to0$ converge weakly to zero, as do their
boundaries; hence $C=0$. Likewise, convergence in $d_V$ implies weak
convergence of cycles, and therefore lower semicontinuity of mass.

Let $(A_j)$ be Cauchy in $\mathcal Z_t$. Pass to a subsequence and
choose fillings $R_j$ of $A_{j+1}-A_j$ with
$\sum_j\M(R_j)<\infty$. The partial sums
$Q_j=\sum_{i=1}^j R_i$ converge in mass to a current $Q$, and
\[
 \partial Q_j=A_{j+1}-A_1,\qquad
 \sup_j\M(\partial Q_j)\leq2t.
\]
By the closure theorem \cite[Theorem~8.5]{AK00}, $Q$ is integral.
Thus $A=A_1+\partial Q$ is an integral cycle,
$\M(A)\leq t$ by lower semicontinuity, and
\[
 d_V(A_{j+1},A)\leq\M(Q-Q_j)\longrightarrow0.
\]
Since the original sequence is Cauchy, it converges to $A$ as well.
\end{proof}

\begin{lemma} \label{lem:equiv}
For any $n\geq 1$,  Theorem \ref{thm:main} is equivalent to Proposition \ref{prop:variational}. 
\end{lemma}

\begin{proof}
	If  Theorem \ref{thm:main} holds true, then we deduce    Proposition \ref{prop:variational} by putting $B=0$.
	
Assume, on the other hand, that Theorem \ref{thm:main} fails 
(for given $n\geq 1$ and $Y$).  Hence, there is a cycle $C\in \I_n(Y)$ and some $c>c_n$ with 
\[
V(C)> c\cdot \M(C)^{\frac {n+1} n}.
\]
Set  $t=\M(C)$.   The function $B\to f(B)= c\cdot \M(B) ^{\frac {n+1} n}$ is a lower semicontinuous function with minimum $0$ on the complete metric space
$(\mathcal Z_t,d_V)$ from Lemma~\ref{lem:filling-metric}. 

Choose $\rho$ with $V(C) >\rho>f(C)$.
 Ekeland's variational principle \cite{Ek74} at $C$
(using parameters $\varepsilon=\lambda=\rho$) provides  an integral cycle
$A\in\mathcal Z_t$ such that
\[
 d_V(A,C)\le\rho,
 \qquad
 f(A)\le f(B)+d_V(A,B)
 \quad\text{for every }B\in\mathcal Z_t.
\]
Since $d_V(C,0)=V(C)>\rho$, we have $A\ne0$.
Furthermore, if $B\in \I_n(Y)$ is a cycle with
$\M(B)>t$, then $f(B)>f(A)$. Consequently,
\[
 c\cdot \bigl(\M(A)^{\frac {n+1}n}-\M(B)^{\frac{n+1} n}\bigr)\le V_Y(A-B) 
\]
for \emph{every} integral $n$-cycle $B\in\I_n(Y)$.
Hence, Proposition~\ref{prop:variational} is violated  as well.
\end{proof}

\section{The inductive argument and setting}\label{sec:compact-support}

We prove Proposition~\ref{prop:variational} and Theorem~\ref{thm:main}
by induction on the cycle dimension $n$.   In the case $n=1$, Theorem~\ref{thm:main} is classical: it follows from
Reshetnyak's majorization theorem and the decomposition of integral
one-cycles into Lipschitz loops; see Subsection \ref{subsec:123} below for details. Applying Lemma \ref{lem:equiv} we deduce Proposition~\ref{prop:variational} for $n=1$.

Let now $n\geq 2$ be fixed and assume that  Proposition~\ref{prop:variational} and Theorem~\ref{thm:main}   hold in all smaller
positive dimensions.  
 Due to  Lemma \ref{lem:equiv}, we only need to verify 
Proposition~\ref{prop:variational} for this dimension $n$.
  We assume  it fails and  find a complete
$\CAT(0)$ space $Y$, a nonzero integral $n$-cycle $A$, and a number
$c>c_n$ such that, 
\begin{equation}\label{eq:variational-compact}
 c\bigl(\M(A) ^{\frac {n+1} {n}}-\M(B)^{\frac {n+1} {n}}\bigr)\leq V_Y(A-B)
 \qquad\text{for every integral $n$-cycle $B$ in $Y$.}
\end{equation}

Denote the volume  of the Euclidean unit
$n$-sphere by
$$s_n=(n+1)\cdot \omega _{n+1}=\frac 1 {c_n ^n \cdot (n+1)^n}\,.$$ 
Rescaling the metric on $Y$ leaves \eqref{eq:variational-compact}
unchanged.  We rescale so that
\begin{equation}\label{eq:mass-normalization}
	\frac 1 {c ^n \cdot (n+1)^n} <m:=\M(A)<s_n.
\end{equation}
We fix this space $Y$, the cycle $A$ and its mass measure $\mu:=\|A\|$  until the end of Section \ref{sec:concentration-completion}. 


After this normalization we linearize \eqref{eq:variational-compact} in the following way:
\begin{lemma}\label{lem:linearized-variation}
There exists $\delta>0$ such that every integral $n$-cycle $B$ in $Y$
satisfies
\begin{equation}\label{eq:linearized-variation}
  \M(A)-\M(B)\leq nV_Y(A-B)
  \qquad\text{whenever }V_Y(A-B)<\delta.
\end{equation}
\end{lemma}

\begin{proof}
Put $m_0=((n+1)c)^{-n}$ and $\delta=c(m^{\frac {n+1} n}-m_0^{\frac {n+1} n})>0$. Let $v:=V_Y(A-B)<\delta$.

There is nothing to prove if $\M(B)\geq m$.  Otherwise, \eqref{eq:variational-compact} gives
\[
  \M(B)^{\frac {n+1} n}\geq m^{\frac {n+1} n}-\frac v c>m_0^{\frac {n+1} n}.
\]
Thus $\M(B)>m_0$, and
\begin{equation} \label{eq:mamb}
  v\geq c\bigl(\M(A)^{\frac {n+1} n}-\M(B)^{\frac {n+1} n}\bigr)
    =c\frac {n+1} n \int_{\M(B)}^{\M(A)} t^{1/n}\,dt
    \geq c\frac {n+1} n m_0^{1/n}\bigl(\M(A)-\M(B)\bigr).
\end{equation}
Since $c\cdot (n+1) \cdot m_0^{\frac 1 n} =1$, inequality  \eqref{eq:mamb}  proves \eqref{eq:linearized-variation}.
\end{proof}


\begin{proposition}
\label{prop:compact-support}
For every $\eta >0$ there exists $a_\eta >0$ with the following properties.
\begin{enumerate}
	\item For every
	Borel set $E\subset Y$ with 
	$\M(A\llcorner E)\leq a_\eta$,   we have
	\begin{equation}\label{eq:uniform-small-isoperimetry}
		\M\bigl(\partial(A\llcorner E)\bigr)
		\geq(1-\eta)n\omega_n^{1/n}\|A\|(E)^{(n-1)/n}.
	\end{equation}
	\item For every $x\in \spt A$ and every $r<a_{\eta}$ we have
	\begin{equation}\label{eq:uniform-growth}
		\|A\|(B(x,r))\geq (1-\eta) \cdot \omega_n \cdot  r^n.
	\end{equation}
\end{enumerate}

$\spt A$ is compact.
\end{proposition}

\begin{proof}
It suffices to consider $0<\eta<1$. Put for this proof
\[
 q=\frac{n+1}{n},\qquad
 \sigma=(1-\eta)^{1/n},\qquad C_0=nD_n2^q.
\]
Choose $b>0$ sufficiently small that
\[
 D_n(2b)^q<\delta,\qquad
 C_0b^{1/n}\leq1-\sigma^{n/(n-1)},
\]
where $\delta$ is the constant in
\eqref{eq:linearized-variation}.

Let $U=A\llcorner E$ with $a:=\M(U)\leq b$. If $\M(\partial U)=\infty $ or if $a=0$, \eqref{eq:uniform-small-isoperimetry} is obvious. Otherwise, the boundary-rectifiability theorem \cite[Theorem~8.6]{AK00} implies  $U\in \I_n (Y)$ and we  
find  a mass-minimizing filling $R\in\I_n(Y)$ of $\partial U$
by Proposition~\ref{prop:wenger-fillings}. Put
$d=\M(R)=V_Y(\partial U)\leq a$.
The current $B=A-U+R$ is an integral cycle, and
\[
 \M(B)\leq m-a+d,\qquad
 V_Y(A-B)=V_Y(U-R)\leq D_n(a+d)^q
                       \leq D_n(2a)^q<\delta.
\]
Consequently, \eqref{eq:linearized-variation} gives
\[
 a-d\leq m-\M(B)\leq nV_Y(A-B)\leq C_0a^{1+1/n}.
\]
On the other hand, the inductive hypothesis in dimension $n-1$
yields
\[
 d=V_Y(\partial U)
 \leq\left(\frac{\M(\partial U)}{n\omega_n^{1/n}}\right)^{n/(n-1)}.
\]
Combining these inequalities and using our choice of $b$, we obtain
\[
 \M(\partial U)
 \geq n\omega_n^{1/n}a^{(n-1)/n}
          (1-C_0a^{1/n})^{(n-1)/n}
 \geq\sigma n\omega_n^{1/n}a^{(n-1)/n}.
\]
Since $\sigma\geq1-\eta$, this proves (1) with any threshold
$a_\eta\leq b$. 

To prove (2), fix $x\in\spt A$, and put
$\beta(r)=\|A\|(B(x,r))$.  By the slicing theorem
\cite[Theorems~5.6--5.7]{AK00}, for almost every $r>0$ we have
\[
 A\llcorner B(x,r)\in\I_n(Y),\qquad
 \M\bigl(\partial(A\llcorner B(x,r))\bigr)\leq\beta'(r).
\]
Thus, whenever $\beta(r)<b$, the estimate established in (1) gives, for almost every such $r$,
\[
 \beta(r)\leq C_\eta\,[\beta'(r)]^{n/(n-1)},\qquad
 C_\eta=(\sigma n\omega_n^{1/n})^{-n/(n-1)}.
\]
Since $x\in\spt A$, we have $\beta(s)>0$ for every $s>0$.
The differential-inequality argument of \cite[Lemma~3.3]{W05},
applied on $[s,r]$ and then letting $s\downarrow0$, therefore yields
\[
 \beta(r)\geq\frac{r^n}{n^n C_\eta^{n-1}}
 =\sigma^n\omega_n r^n=(1-\eta)\omega_n r^n
\]
whenever $\beta(r)<b$.  Decreasing $a_{\eta}$ if necessary we obtain (2).

Applying (2) with $\eta=1/2$ and using
$\mu(Y)<\infty$ excludes infinitely many uniformly separated
points in $\spt A$. Thus $\spt A$ is totally bounded, hence compact.
\end{proof}

\section{Infinitesimal estimates for radial deformations}
\label{sec:radial}

We record the main deformation estimate.  Let $0\ne C\in\I_n(Y)$ have
compact support $K$, and put $\mu=\|C\|$.  Let
$g:K\to[0,\infty)$ be Lipschitz.  Fix $y\in Y$ and let $r=d_y$
denote the distance function to $y$.

For every $z\in Y$, denote by $\gamma^z:[0,1]\to Y$ the
constant-speed geodesic from $y$ to $z$.  Consider the radial
contraction $F=F^{y,g}:[0,1]\times K\to Y$, defined by
\[
 F(t,x):=\gamma^x\bigl(\max\{1-tg(x),0\}\bigr).
\]
The map $F$ is Lipschitz by the CAT(0) property.
For $h\in[0,1]$, set $F_h(x)=F(h,x)$. 
We  endow $[0,1]\times Y$ with the product metric and consider 
  integral currents
\[
 B_h=(F_h)_\#C\in\I_n(Y),\qquad
 D_h=F_\#(\llbracket0,h\rrbracket\times C)\in\I_{n+1}(Y).
\]
If $C$ is a cycle, then $B_h$ is a cycle as well, and the homotopy
formula gives
\[
 \partial(\llbracket0,h\rrbracket\times C)
 =\llbracket h\rrbracket\times C-\llbracket0\rrbracket\times C,
 \qquad
 \partial D_h=B_h-C.
\]

\begin{proposition}
\label{prop:radial-current-bounds}
In the notation above, as $h\downarrow0$,
\begin{align}
 \M(B_h)&\leq\M(C)
  -h\int_K\bigl(ng+r\langle\nabla^C r,\nabla^C g\rangle\bigr)
       \,d\mu+O(h^2),
       \label{eq:radial-first-variation}\\
 \M(D_h)&\leq h\int_K rg\sqrt{1-|\nabla^C r|^2}\,d\mu.
       \label{eq:radial-swept-first-order}
\end{align}
\end{proposition}

\begin{proof}
Fix $0<h\leq1$ with $h\|g\|_\infty<1$, and put
$\widehat C_h=\llbracket0,h\rrbracket\times C$.
Then
\[
 \|\widehat C_h\|
 = (\LL^1\llcorner(0,h))\otimes\mu,
 \qquad T_{(t,x)}\widehat C_h=\R\oplus T_xC
 \quad\text{almost everywhere}.
\]
Let $J(x,h)$ and $\widetilde J(t,x)$ be the Jacobians of
$d^C(F_h)_x$ and $d^{\widehat C_h}F_{(t,x)}$, respectively.
By the area formula \cite{AK00,K94} and \eqref{eq:mass-formula},
\[
 \M(B_h)\leq\int_K J(x,h)\,d\mu,
 \qquad
 \M(D_h)\leq\int_0^h\int_K\widetilde J(t,x)\,d\mu\,dt.
\]

For $t,t'\in[0,h]$ and $x,x'\in K$, set
$\alpha=1-tg(x)$ and $\alpha'=1-t'g(x')$.
Comparison with the Euclidean triangle for $y,x,x'$ gives
\begin{equation}\label{eq:radial-comparison}
 d(F(t,x),F(t',x'))^2
 \leq\alpha\alpha'd(x,x')^2
     +(\alpha-\alpha')\bigl(\alpha r(x)^2-\alpha'r(x')^2\bigr).
\end{equation}
At $\mu$-almost every $x$, write
$a=\nabla^Cr(x)$, $b=\nabla^Cg(x)$ and $\alpha_t=1-tg(x)$,
suppressing the argument $x$ in $r$ and $g$ below.
To differentiate \eqref{eq:radial-comparison}, fix a density point
$z$ in a rectifiable chart $\phi$, with $x=\phi(z)$, and choose
$w$ with $d\phi_z(w)=v$. For $0<t<h$, take
$t'=t+\varepsilon\tau$ and chart points $z'$ satisfying
$z'-z=\varepsilon w+o(\varepsilon)$ as $\varepsilon\downarrow0$.
Writing $x'=\phi(z')$, we have
\[
\begin{aligned}
	d(x,x')&=\varepsilon|v|+o(\varepsilon),&
	r(x')&=r+\varepsilon\langle a,v\rangle+o(\varepsilon),&
	\alpha'&=\alpha_t-\varepsilon (t\langle b,v\rangle+g\tau)+o(\varepsilon).
\end{aligned}
\]
Let $\eta=t\langle b,v\rangle+g\tau$ be the first-order coefficient
in the last expansion. Inserting these expansions into
\eqref{eq:radial-comparison}, dividing by $\varepsilon^2$
and passing to the limit yields for almost every $(t,x)$
\begin{equation}  \label{eq:longeq}
|d^{\widehat C_h}F_{(t,x)}(\tau,v)|^2
\leq\alpha_t^2|v|^2-2\alpha_t r\eta\langle a,v\rangle+r^2\eta^2.
\end{equation}
We define a comparison linear map 
$L_t:\R \oplus T_xC \to \R \oplus T_xC$ by
\[
 L_t(\tau,v)=
 \bigl(\ r\eta\sqrt{1-|a|^2}, \alpha_t v-r\eta a\bigr).
\]
Then $|L_t(\tau,v)|^2$ equals the right hand side of \eqref{eq:longeq}. 
The same argument with $t=t'=h$ gives
$|d^C(F_h)_x(v)|\leq|L_h(0,v)|$ for $\mu$-almost every $x$.
These differentials have Euclidean images by
Proposition~\ref{prop:tangential-differentials}; hence, with
$M_h=L_h|_{\{0\} \oplus T_xC}$, their Jacobians satisfy
\begin{equation}\label{eq:newcomp}
 J(x,h)\leq J_n(M_h),
 \qquad \widetilde J(t,x)\leq|\det L_t|.
\end{equation}

Multilinearity of the determinant gives
\begin{equation}\label{eq:pointwise}
 \widetilde J(t,x)\leq|\det L_t|
 =\alpha_t^n rg\sqrt{1-|a|^2}
 \leq rg\sqrt{1-|a|^2}.
\end{equation}
In order to estimate the Jacobian of $M_h$, we note that
\[
\begin{aligned}
 M_0(v)&=(0,v),\\
 M'_0(v)&=\bigl(r\langle b,v\rangle\sqrt{1-|a|^2},
                 -gv-r\langle b,v\rangle a\bigr).
\end{aligned}
\]
Thus the first-order expansion of
$J_n(M_h)$ gives
\begin{equation}\label{eq:radial-integrated-jacobian}
 J(x,h)\leq J_n(M_h)
 =1-h\bigl(ng+r\langle a,b\rangle\bigr)+O(h^2).
\end{equation}
The remainder is uniform in $x$, since $|a|\leq1$,
$|b|\leq\Lip(g)$, and $r,g$ are bounded on $K$.
Integrating \eqref{eq:pointwise} and
\eqref{eq:radial-integrated-jacobian} proves  \eqref{eq:radial-first-variation}   and \eqref{eq:radial-swept-first-order}.
\end{proof}

The "extremality"  of  $A$, expressed by \eqref{eq:linearized-variation} together with Proposition \ref{prop:radial-current-bounds} implies:

\begin{corollary}
\label{cor:radial-variation}
For the cycle $A$ fixed in Section~\ref{sec:compact-support},
write $\mu=\|A\|$ and $r(x)=d(y,x)$, where $y\in Y$ is arbitrary.
Then every nonnegative Lipschitz function $g$ on $\spt A$ satisfies
\begin{equation}\label{eq:radial-variational-inequality}
 \int\bigl(ng+r\langle\nabla^A r,\nabla^A g\rangle\bigr)\,d\mu
 \leq n\int rg\sqrt{1-|\nabla^A r|^2}\,d\mu.
\end{equation}
\end{corollary}

\begin{proof}
Apply Proposition~\ref{prop:radial-current-bounds} with $C=A$.
Since $\partial D_h=B_h-A$, estimate
\eqref{eq:radial-swept-first-order} gives
$V_Y(A-B_h)\leq\M(D_h)\to0$ as $h\downarrow0$.
Thus \eqref{eq:linearized-variation} applies for all sufficiently
small $h$ and yields
\[
 \M(A)-\M(B_h)\leq n\M(D_h).
\]
Combine this with \eqref{eq:radial-first-variation} and
\eqref{eq:radial-swept-first-order}, divide by $h$, and let
$h\downarrow0$.
\end{proof}

\begin{remark}\label{rem:convex-flow-variation}
The radial inequality has an analogue for general $\lambda$-convex potentials.
Let $\varphi:Y\to\R$ be a locally Lipschitz $\lambda$-convex function,
let $\Phi_t$ denote its descending gradient flow, and write
$\sigma=|\partial^-\varphi|$ for its descending slope.
Using the deformation $F(t,x)=\Phi_{tg(x)}(x)$ in place of the radial
contraction and relying on the 
 unequal-time comparison for these flows established in
\cite[Lemma~3.4]{GN21}, one obtains the following conclusion: 
for  every nonnegative Lipschitz function $g$ on
$K=\spt A$,
\[
 \int_K\bigl(n\lambda g+
      \langle\nabla^A\varphi,\nabla^A g\rangle\bigr)\,d\mu
 \leq n\int_K g\sqrt{\sigma^2-|\nabla^A\varphi|^2}\,d\mu.
\]
Taking $\varphi=\frac12 d_y^2$ gives $\lambda=1$, $\sigma=r$ and
$\nabla^A\varphi=r\nabla^A r$, recovering
Corollary~\ref{cor:radial-variation}.
\end{remark}

\section{Dimension 2}\label{sec:dimension-two}

We prove Proposition~\ref{prop:variational} for $n=2$, following
\cite[Proposition~3.2]{OAI26}. We retain the setting of
Section~\ref{sec:compact-support} and put $K=\spt A$ and $\mu=\|A\|$.
In particular, \eqref{eq:mass-normalization} gives
\[
  0<m=\mu(K)<s_2=4\pi.
\]
Fix any point $y\in K$. By Proposition~\ref{prop:compact-support},
\begin{equation}\label{eq:dimension-two-density}
  \liminf_{\varepsilon\downarrow0}
  \frac{\mu(B(y,\varepsilon))}{\varepsilon^2}\geq\omega_2=\pi.
\end{equation}
Put $r=d_y$. For $0<\varepsilon<1$ with
$\mu(\{r=\varepsilon\})=0$, consider the nonnegative Lipschitz function
\[
  g_\varepsilon(x)=\max\{r(x),\varepsilon\}^{-2},\qquad x\in K.
\]
Almost everywhere on $K$, the chain rule gives
\[
  \nabla^A g_\varepsilon=
  \begin{cases}
    0,&r<\varepsilon,\\
    -2r^{-3}\nabla^A r,&r>\varepsilon.
  \end{cases}
\]

Applying Corollary~\ref{cor:radial-variation} to $g_\varepsilon$,
splitting the integrals over the  ball $B(y,\varepsilon)$  and its complement,  and dividing by 2, we obtain
\[
  \frac{\mu(B(y,\varepsilon))}{\varepsilon^2}
  +\int_{\{r>\varepsilon\}}
        \frac{1-|\nabla^A r|^2}{r^2}\,d\mu
  \leq
  \int_{\{r<\varepsilon\}}
        \frac r{\varepsilon^2}\sqrt{1-|\nabla^A r|^2}\,d\mu
  +\int_{\{r>\varepsilon\}}
        \frac1r\sqrt{1-|\nabla^A r|^2}\,d\mu.
\]
The first integrand on the right is bounded by $1/\varepsilon$.  We deduce 
\[
\left(\frac1{\varepsilon^2}-\frac1\varepsilon\right)
\mu(B(y,\varepsilon))
\leq\int_{\{r>\varepsilon\}}(f-f^2)\,d\mu,
\]
where we put
\[
  f=\frac1r\sqrt{1-|\nabla^A r|^2}
  \qquad\text{on }\{r>\varepsilon\}.
\]
Since $f-f^2\leq\frac14$, we conclude that
\begin{equation}\label{eq:dimension-two-key-estimate}
  (1-\varepsilon)\frac{\mu(B(y,\varepsilon))}{\varepsilon^2}
  \leq\frac14\mu(\{r>\varepsilon\})\leq\frac m4.
\end{equation}
Only countably many level sets of $r$ have positive $\mu$-measure.
Hence we can let $\varepsilon\downarrow0$ through radii satisfying
$\mu(\{r=\varepsilon\})=0$. By \eqref{eq:dimension-two-density},
we obtain $\pi\leq m/4$, contradicting $m<4\pi$.

This proves Proposition~\ref{prop:variational}, and therefore
Theorem~\ref{thm:main} for $n=2$.

\section{Sobolev spaces on integral currents}\label{sec:sobolev-spaces}

Let $0\ne C\in\I_n(Y)$ be an arbitrary integral current with compact
support $K=\spt C$, and put $\mu=\|C\|$.  We do not assume that
$\partial C=0$.  For Lipschitz functions $u$ on $K$, let $d^Cu$
denote the $\mu$-almost everywhere defined tangential differential
introduced in Section~\ref{sec:prelim}.  Its $L^2$-closure will be denoted by $D^C$.
For $u\in\Lip(K)$, set
\[
 \|u\|_{\mathcal H_C}^2=\|u\|_2^2+\|d^Cu\|_2^2,\qquad
 \mathcal H_C=\overline{\Lip(K)}^{\,\|\cdot\|_{\mathcal H_C}}.
\]
All function norms in this section are taken with respect to $\mu$.

The general theory \cite[Corollary~4.3 and Theorems~5.3, 8.1]{P15}
shows that $d^C$ is closable in $L^2(\mu)$ and that
$\mathcal H_C$ identifies with the Cheeger--Sobolev space
$W^{1,2}(K,d,\mu)$.  In the normalization in which the Cheeger
energy $\operatorname{Ch}_2$ contains a factor $1/2$,
every $u\in\mathcal H_C$ satisfies
\[
 \|u\|_{\mathcal H_C}^2=\|u\|_2^2+\|D^Cu\|_2^2,\qquad
 \|D^Cu\|_2^2=2\operatorname{Ch}_2(u),\qquad
 |D^Cu|=g_u\quad\mu\text{-a.e.},
\]
where $g_u$ is the minimal weak upper gradient.  The equivalence
with Newtonian Sobolev spaces and Lipschitz density are covered by
\cite{AGS13}; in particular, $\mathcal H_C=N^{1,2}(K,d,\mu)$,
with equality of norms.  If the approximate tangent norms of $C$
are induced by inner products $\mu$-almost everywhere, then
$\mathcal H_C$ is Hilbertian \cite[Theorem~8.1]{P15}.
This condition holds in our $\CAT(0)$ setting by
Theorem~\ref{thm:charts}; see also \cite[Theorem~1.1]{DGPS21}.

We use the standard Sobolev calculus \cite{AGS13,HKST15}.  For
$u_1,\ldots,u_k\in\mathcal H_C$ and $\Phi\in C^1(\R^k)$ with
bounded first derivatives,
\begin{equation}\label{eq:closed-chain-rule}
 \Phi(u_1,\ldots,u_k)\in\mathcal H_C,\qquad
 D^C\Phi(u_1,\ldots,u_k)
 =\sum_{j=1}^k\partial_j\Phi(u_1,\ldots,u_k)D^Cu_j.
\end{equation}
The differential of $u\in \mathcal H_C$ vanishes $\mu$-almost everywhere on level sets
$\{u=t\}$.  The scalar chain rule extends to Lipschitz piecewise-$C^1$
functions.  In
particular, $u\mapsto u_+$ is continuous on $\mathcal H_C$, and
nonnegative elements admit nonnegative Lipschitz approximations.
We refer to these properties and the closedness of $D^C$
collectively as the \emph{Sobolev calculus}.

\begin{proposition}\label{prop:intrinsic-rellich}
For every $1\leq q<2$, the canonical embedding
$\mathcal H_C\hookrightarrow L^q(\mu)$ is compact.
\end{proposition}

\begin{proof}
It suffices to prove compactness into $L^1(\mu)$; the uniform
$L^2$-bound and H\"older interpolation
\cite[Section~6.5]{Folland} then give compactness into $L^q(\mu)$
for every $1\leq q<2$.

Rescaling the metric, we  assume $\mu (K)=1$.
 By Lipschitz density, it suffices  to show that every sequence \((u_j)\subset\operatorname{Lip}(K)\) with \(\|u_j\|_{\mathcal H_C}\le1\) admits an \(L^1(\mu)\)-convergent subsequence.

 Consider, for any $u\in \mathcal H_C$ and $M>0$, the truncation
$F_M(u)=\max\{-M,\min\{u,M\}\}$.  Then $F_M(u)\in \mathcal H_C$ and $\|F_M(u)\|_{\mathcal H_C } \leq \|u\| _{\mathcal H_C}$.  Moreover,
\[\|u-F_M(u)\|_1  \leq \int_ {|u|>M} ( |u|-M) \, d\mu  \leq \frac 1 M \|u\|_2^2\,.\]  
Thus,  if for every $M>0$  the sequence $F_M(u_j)$ is precompact in $L^1(\mu)$  then, a diagonal argument provides an $L^1$-convergent subsequence of the original sequence $u_j$.    Thus, without loss of generality, we may assume that our Lipschitz functions $u_j$ are uniformly bounded. Upon rescaling we may assume $\|u_j\|_{\infty}\leq 1$ for all $j$.

We claim that it is sufficient to verify that $\mathbf 1_Z \cdot u_j $ 
has a convergent subsequence, 
 for any image $Z$ of a bilipschitz chart $\phi:E\to K$, with $E\subset \R^n$  compact. Indeed, up to a $\mu$-negligible set we have a disjoint subdivision of $K$  in such images $Z_i$.
 Moreover, for any $j$:
 \[\|u_j -(\sum _{i=1} ^k  \mathbf 1_{Z_i}\cdot u_j)\|_1 \leq \|u_j\|_{\infty} \cdot \mu (K\setminus \cup _{i=1} ^{k}  Z_i) \leq  \mu (K\setminus \cup _{i=1} ^{k}  Z_i).\]
 Thus,  knowing that the restriction of $(u_j)$ to any $Z_i$ has a convergent subsequence would produce a convergent subsequence in $L^1(\mu)$  by a diagonal argument.
 
 We  fix a bilipschitz chart  $\phi :E\to Z\subset K$ as above and need to verify that  $\mathbf 1_Z \cdot u_j $  is precompact in $L^1(\mu)$.  Moreover, by the same decomposition argument as above, we may assume that the multiplicity $\theta$ of $C$ is constant on $Z$.
Choose the orientation of the chart so that $\theta>0$.
Extend $\phi^{-1}:Z\to E$ to a Lipschitz map $F:K\to\R^n$.

For any fixed $\chi\in\Lip(K)$ with $0\le\chi\le1$, consider the normal current
\[
 S_{j,\chi}:=F_\#(C\llcorner(\chi u_j)).
\]
The product rule and our normalizations give
\[
 \M(C\llcorner(\chi u_j))
 +\M(\partial(C\llcorner(\chi u_j)))
 \le \M(\partial C)+2+\Lip(\chi).
\]
Thus the normal $n$-currents $S_{j,\chi}$ have uniformly bounded
mass and boundary mass, and are supported in the compact set $F(K)$.
Since they are top-dimensional currents in $\R^n$, the 
$BV$ compactness theorem \cite[Theorem~3.23]{AFP00} implies that
$(S_{j,\chi})$ is precompact in mass.

Now put $\chi_k=\max\{1-kd(\cdot,Z), 0\}$ and 
\[
 S_{j,Z}:=F_\#(C\llcorner(\mathbf1_Zu_j)).
\]
Since $Z$ is compact and $|u_j|\le1$,
\[
 \sup_j\M(S_{j,\chi_k}-S_{j,Z})
 \le\Lip(F)^n\|\chi_k-\mathbf1_Z\|_1\longrightarrow0.
\]
Hence $(S_{j,Z})$ is precompact in mass as well.

By \eqref{eq:representation}, the constant multiplicity $\theta>0$,
and $F\circ\phi=\mathrm{id}_E$,
\[
 \M(S_{j,Z}-S_{\ell,Z})
 =\theta\int_E|u_j\circ\phi-u_\ell\circ\phi|\,dz.
\]
The area formula therefore gives
\[
 \|(u_j-u_\ell)\mathbf1_Z\|_1
 \le\Lip(\phi)^n\M(S_{j,Z}-S_{\ell,Z}).
\]
This proves the required precompactness and completes the proof.
\end{proof}

\begin{corollary}  \label{cor:radial-sobolev}
Under the assumptions of Corollary~\ref{cor:radial-variation}, every nonnegative 
 $g\in\mathcal H_A$ satisfies
\begin{equation}\label{eq:radial-variational-sobolev}
 \int\bigl(ng+r\langle Dr,Dg\rangle\bigr)\,d\mu
 \leq n\int rg\sqrt{1-|Dr|^2}\,d\mu;
\end{equation}
\begin{equation}\label{eq:radial-quadratic-sobolev}
 \int\bigl(ng|Dr|^2+r\langle Dr,Dg\rangle\bigr)\,d\mu
 \leq\frac n4\int r^2g\,d\mu.
\end{equation}
\end{corollary}

\begin{proof}
Approximating any $g$ by nonnegative Lipschitz functions in $\mathcal H$ we directly deduce \eqref{eq:radial-variational-sobolev} from
 Corollary~\ref{cor:radial-variation}. The inequality \eqref{eq:radial-quadratic-sobolev} follows from \eqref{eq:radial-variational-sobolev}  and  the pointwise identity
\[
 1-r\sqrt{1-|Dr|^2}-|Dr|^2+\frac{r^2}{4}
 =\left(\sqrt{1-|Dr|^2}-\frac r2\right)^2\geq0 .
\]
\end{proof}

\begin{remark}
	Proposition~\ref{prop:intrinsic-rellich} uses neither $\CAT(0)$
	geometry nor the cycle condition.  It holds for every compactly
	supported integral current in an arbitrary complete metric space; Euclidean tangent norms
	are needed for Hilbertianity, but not for this compactness result.
	Alternatively, its conclusion might be obtained by combining the
	local Poincar\'e estimates of \cite[Theorem~C.1]{P15} with the
	abstract compactness criterion of \cite[Theorem~2.3]{BKa22}.
\end{remark}

\section{Small-set isoperimetry and the Sobolev inequality}
\label{sec:intrinsic-sobolev}

We now return to the cycle $A\in\I_n(Y)$ fixed in
Section~\ref{sec:compact-support} and assume $n\geq3$. Its support
$K=\spt A$ is compact by Proposition~\ref{prop:compact-support}.
As before, put $\mu=\|A\|$, $m=\mu(K)=\M(A)$, and
$\mathcal H=\mathcal H_A$, $D=D^A$, using the notation of
Section~\ref{sec:sobolev-spaces}. 
We write the critical Sobolev exponent and the Euclidean constants as
\begin{equation}\label{eq:intrinsic-sobolev-constants}
  p=\frac{2n}{n-2},\qquad k_n=ns_n^{2/n}.
\end{equation}
Thus, the optimal Sobolev inequality  for $f\in W^{1,2}(\R^n)$ due to Aubin--Talenti \cite{Aubin76}, \cite{Talenti} is
\begin{equation} \label{eq:talenti}
k_n\|f\|_{L^p(\R^n)}^2\leq\frac4{n-2}\|\nabla f\|_{L^2(\R^n)}^2.
\end{equation}

We derive the Sobolev inequality on $(K,d,\mu)$ from the metric-measure generalization of the classical Pólya--Szegő inequality \cite{PS51}, established by Nobili--Violo \cite{NV25}.
 In order to use this version, we recall the notion of perimeter in metric measure spaces and compare it to the boundary mass.  For a Borel subset $E\subset K$ its perimeter is   \cite[Definition~2.2]{NV25} 
\begin{equation} \label{eq:perdef}
	\operatorname{Per}_\mu(E)
	=
	\inf_{f_j\to\mathbf1_E}
	\liminf_{j\to\infty}
	\int_K\operatorname{lip}f_j\,d\mu,
\end{equation}
where the infimum is taken over Lipschitz functions  \(f_j:K\to [0,1]\) converging to the characteristic function  \(\mathbf1_E\) in \(L^1(\mu)\).
Here $\operatorname{lip}f_j$ is the pointwise Lipschitz constant, which in our case satisfies  $|D f_j|\leq \operatorname{lip}f_j$ almost everywhere on $K$.

We can use a diagonal argument and  find a sequence of Lipschitz functions $f_j$ as above,  such that $\operatorname{Per}_\mu(E)
=\lim_{j\to\infty}
\int_K\operatorname{lip}f_j\,d\mu$. The restrictions $A\llcorner f_j$ are normal currents converging to $A\llcorner E$.
Since $\partial A =0$, the boundary formula for metric currents implies 
\[\M(\partial (A\llcorner f_j))= \M (A\llcorner df_j) \leq   \int _K |D f_j| \, d\mu  \leq \int _K \operatorname{lip}f_j \, d\mu.  
\]
From  the lower semicontinuity of boundary mass \cite{AK00}, we derive 
 \begin{equation}
\M\bigl(\partial(A\llcorner E)\bigr) \leq \liminf _{j\to \infty} \M(\partial(A\llcorner f_j)) 
\leq\operatorname{Per}_{\mu}(E).
\end{equation}
Combining this comparison with Proposition~\ref{prop:compact-support}, we obtain, for every $\eta\in (0,1)$ and every Borel set $E\subset K$ with $\mu(E)\leq a_\eta$, the estimate
\begin{equation}\label{eq:intrinsic-bv-isoperimetry}
	\operatorname{Per}_{\mu}(E)
	\geq(1-\eta)n\omega_n^{1/n}\mu(E)^{(n-1)/n}.
\end{equation}

\begin{proposition}
\label{prop:intrinsic-sobolev}
For any $\eta\in(0,1)$ and any nonnegative Lipschitz function
$u:K\to\R$ with $\mu\{u>0\}\leq a_\eta$,
\begin{equation}\label{eq:small-support-sobolev}
 (1-\eta)^2k_n\|u\|_p^2\leq\frac4{n-2}\|Du\|_2^2.
\end{equation}
\end{proposition}

\begin{proof}
We apply \cite[Theorem~1.1]{NV25}  with $\Omega=X=K$, $\mathsf C=1-\eta$, with exponent $2$ 
and the weighted interval $I=(0,\infty)$ with measure
$
d\omega(r)=g(r)\,dr$, where 
$g(r)=n\omega_n r^{n-1}$.
The corresponding isoperimetric profile in the
notation of \cite{NV25} is
\[
\mathcal I_g^\flat(s)
=n\omega_n^{1/n}s^{(n-1)/n}.
\]

For $\delta>0$ put $u_\delta=(u-\delta)_+$.
Since a sufficiently small tubular neighborhood of  $\spt u_\delta$ is contained in $\{u>0\}$, the isoperimetric hypothesis
(1.1) of \cite{NV25} follows from
\eqref{eq:intrinsic-bv-isoperimetry}.

Put $R=(m/\omega_n)^{1/n}$ and let
$u_\delta^*:(0,R)\to[0,\infty)$ be the
decreasing rearrangement of $u_\delta$
with respect to $\omega$.
Inequality (1.2) of \cite{NV25} yields
\[
\int_K|Du_\delta|^2\,d\mu
\geq (1-\eta)^2
\int_0^R |(u_\delta^*)'(r)|^2\,d\omega(r).
\]
Here we use the identification of the minimal
weak upper gradient with $|Du_\delta|$
from Section~\ref{sec:sobolev-spaces}.

Define $w_{\delta} :\R^n \setminus \{0 \} \to \R$ by setting $w_\delta(x)=u_\delta^*(|x|)$ for $|x|<R$ and $w_{\delta} (x)=0$ for $|x|\geq R$.
Since necessarily $a_\eta<m$, the rearrangement
$u_\delta^*$ vanishes near $R$.
Therefore, $w_{\delta} \in W^{1,2}(\R^n \setminus \{0\})$.  Since $n\geq 3$, $w_{\delta} \in W^{1,2} (\R^n)$   and 
\[
\int_{\R^n}|\nabla w_\delta|^2\,d\mathcal L^n
=\int_0^R |(u_\delta^*)'(r)|^2\,d\omega(r)
\leq(1-\eta)^{-2}\|Du_\delta\|_2^2.
\]
Since $u^{\ast}_{\delta}$ is the rearrangement of $u_{\delta}$ and by the definition of    $w_{\delta}$, we have
\[\|u_{\delta}\|_{L^p(\mu)}^2=\|u^{\ast}  _{\delta}\|_{L^p(\omega)}^2  = \|w_{\delta}\|_{L^p(\R^n)}^2\,.\]
Hence, the Euclidean Sobolev inequality \eqref{eq:talenti}  implies
\[
\begin{aligned}
	(1-\eta)^2 k_n\|u_\delta\|_p^2
	=(1-\eta)^2 k_n
	\|w_\delta\|_{L^p(\mathbb R^n)}^2
	&\leq \frac{4(1-\eta)^2}{n-2}
	\|\nabla w_\delta\|_{L^2(\mathbb R^n)}^2\\
	&\leq \frac4{n-2}\|Du_\delta\|_2^2
	\leq \frac4{n-2}\|Du\|_2^2.
\end{aligned}
\]
Letting $\delta\downarrow0$ proves \eqref{eq:small-support-sobolev}.
\end{proof}

\begin{corollary}\label{cor:intrinsic-sobolev}
We have $\mathcal H \subset L^p (\mu)$.  For every $1\leq q<p$ the canonical embedding $\mathcal H\to L^q(\mu)$ is compact.
	
	For every $\varepsilon\in(0,k_n)$ there is some $b_\varepsilon\geq0$ such that for all $u  \in \mathcal H$
	\begin{equation}\label{eq:almost-sharp-sobolev}
		(k_n-\varepsilon)\|u\|_p^2
		\leq\frac4{n-2}\|Du\|_2^2+b_\varepsilon\|u\|_2^2.
	\end{equation}	
\end{corollary}

	\begin{proof}
	We may replace $u$ by $|u|$ and assume $u\geq 0$.	Fix $\eta\in(0,1)$. Since $\mu$ is nonatomic
		and $K$ is compact, we find finitely many
		nonnegative Lipschitz functions $\chi_i$ on $K$
		such that
		\[
		\sum_i\chi_i^2=1,\qquad
		\mu(\spt\chi_i)\leq a_\eta.
		\]
		Set $C_{\eta}=\sum_i \Lip(\chi_i)^2 <\infty$. 
		For every Lipschitz function $u$ we have  $u^2=\sum_i (\chi_i u)^2$ and Proposition~\ref{prop:intrinsic-sobolev} applies to every function $\chi_iu$.  Minkowski's
		inequality in $L^{p/2}(\mu)$ gives us 
		\[(1-\eta)^2k_n\|u\|_p^2
			=(1-\eta)^2k_n
			\left\|\sum_i(\chi_i u)^2\right\|_{p/2}
			\leq(1-\eta)^2k_n\sum_i\|\chi_i u\|_p^2\leq\frac4{n-2}\sum_i\|D(\chi_i u)\|_2^2
	\]	
	 $\sum_i \chi_i ^2=1$ implies $\sum _i\chi_i D{\chi_i} =0$. From  the Leibniz rule we obtain almost everywhere			
	\[\sum_i |D(\chi_i u)|^2=|Du|^2\sum_i |\chi_i|^2+ |u|^2 \sum_i |D\chi_i| ^2+2 \langle \sum \chi_i D\chi _i, uDu \rangle \leq |Du|^2+C_{\eta} u^2
		\]
	Taking these inequalities  together we deduce
		\[
		(1-\eta)^2k_n\|u\|_p^2
		\leq\frac4{n-2}\|Du\|_2^2
		+\frac{4C_\eta}{n-2}\|u\|_2^2.
		\]
		Choosing $\eta$ sufficiently small proves
		\eqref{eq:almost-sharp-sobolev} for Lipschitz
		functions. 
		Since Lipschitz functions are dense in $\mathcal H$, the inequality holds true for all $u\in \mathcal H$ and $\mathcal H$ is contained in $L^p(\mu)$.  As in the proof of  Proposition \ref{prop:intrinsic-rellich},  interpolation implies that the embedding $\mathcal H\to L^q(\mu)$ is compact
		 for all $1\leq q< p$.		
\end{proof}

\section{A normalized minimizer and nonlinear test functions}
\label{sec:sobolev-minimizer}
We apply the abstract Yamabe theorem of Akutagawa--Carron--Mazzeo
\cite[Theorem~4.9(i)]{ACM19} to obtain a bounded minimizer of a
Sobolev quotient. On $\mathcal H$, we use the symmetric bilinear form
$Q(\cdot,\cdot)$ and its associated quadratic form:
\begin{equation}\label{eq:spherical-quotient}
 \begin{aligned}
 Q(u,w)&=\frac4{n-2}\int_K\langle Du,Dw\rangle\,d\mu
          +n\int_Kuw\,d\mu,\\
 Q(u)&=Q(u,u),\qquad
 \Lambda=\inf_{0\ne u\in\mathcal H}\frac{Q(u)}{\|u\|_p^2}.
 \end{aligned}
\end{equation}

The following result uses (for the first time) our normalization $m<s_n$.

\begin{proposition}\label{prop:normalized-minimizer}
We have
$0<\Lambda\leq nm^{2/n}<k_n$.
Moreover, there exists a nonnegative function
$u\in\mathcal H\cap L^\infty(\mu)$ satisfying
$\|u\|_p=1$ and $Q(u)=\Lambda$.
\end{proposition}

\begin{proof}
  By the  Sobolev
inequality \eqref{eq:almost-sharp-sobolev}, the quantity
$\|u\|_p^2$ is bounded by a constant multiple of $Q(u)$;
thus $\Lambda>0$.  Testing the quotient with the constant function
$1$ yields
\[
 \Lambda\leq \frac{nm}{m^{2/p}}
 =nm^{2/n}<ns_n^{2/n}=k_n.
\]

To apply \cite[Theorem~4.9(i)]{ACM19}, consider the quadratic form
\[
 \mathcal E(f,g)=\frac4{n-2}\int_K\langle Df,Dg\rangle\,d\mu,
 \qquad \mathcal D(\mathcal E)=\mathcal H.
\]
This is a regular, strongly local Dirichlet form on $K$.  More explicitly:
\begin{enumerate}
  \item $\mathcal H$ is complete with respect to
  $\|f\|_{\mathcal E_1}^2=\mathcal E(f,f)+\|f\|_2^2$;
  \item $\tilde f=\max\{0,\min\{f,1\}\}$ belongs to $\mathcal H$ and
  $\mathcal E(\tilde f,\tilde f)\leq\mathcal E(f,f)$ (the Markov property);
  \item $\mathcal E(f,g)=0$ if $f$ is constant on a neighborhood
  of  $\spt g$ (strong locality);
  \item $\mathcal H\cap C(K)$ is dense in $\mathcal H$ for
  $\|\cdot\|_{\mathcal E_1}$ and in $C(K)$ for the uniform norm
  (regularity).
\end{enumerate}
Indeed, (1) follows from the equivalence of $\|\cdot\|_{\mathcal E_1}$
with the norm on $\mathcal H$ defined in Section~\ref{sec:sobolev-spaces}.
Items (2) and (3) follow  from locality of $D$.
Finally, Lipschitz functions are dense in $\mathcal H$ for the Sobolev
norm and in $C(K)$ for the uniform norm, proving (4).

Choose $\varepsilon\in(0,k_n-\Lambda)$ and rewrite
\eqref{eq:almost-sharp-sobolev} as
\[
 \|f\|_p^2\leq
 \frac1{k_n-\varepsilon}\mathcal E(f,f)
 +\frac{b_\varepsilon}{k_n-\varepsilon}\|f\|_2^2.
\]
This is the Sobolev inequality of \cite{ACM19}, with
$\nu=n$, $A=(k_n-\varepsilon)^{-1}$,
$B=b_\varepsilon/(k_n-\varepsilon)$, and constant potential $V=n$.
The corresponding Yamabe energy is
$\mathcal E_V(f)=Q(f)$ and its Yamabe constant is
$Y(V)=\Lambda<1/A=k_n-\varepsilon$.
Since $V_-=0\in L^q(\mu)$,
all the hypotheses of \cite[Theorem~4.9(i)]{ACM19} are satisfied.
The theorem provides a nonnegative bounded minimizer $u$ with
$\|u\|_p=1$ and $Q(u)=\Lambda$.
\end{proof}

Let $u$ be the nonnegative normalized minimizer provided by
Proposition~\ref{prop:normalized-minimizer}.
Define
\[
 c=\left(\frac{\Lambda}{n}\right)^{\frac1{p-2}},\qquad v=cu.
\]
The bounded nonnegative function $v$ is still a minimizer, and
$\|v\|_p=c$. In particular,
\begin{equation}\label{eq:normalized-weighted-mass}
 \begin{aligned}
 0<W:=\int_Kv^p\,d\mu
=\|v\|_p^p=c^p
 =\left(\frac{\Lambda}{n}\right)^{\frac{p}{p-2}}
 =\left(\frac{\Lambda}{n}\right)^{\frac n2}
 \leq m<s_n.
 \end{aligned}
\end{equation}
Since $Q(u)=\Lambda$ and $c^{p-2}=\Lambda/n$, the quadratic
homogeneity of $Q$ gives
\begin{equation}
 Q(v)=c^2Q(u)=c^2\Lambda=nc^p=nW.
\end{equation}

We fix this choice of $v$ and derive its Euler--Lagrange equation.
\begin{lemma}\label{lem:normalized-euler-lagrange}
For every $\psi\in\mathcal H$, we have
\begin{equation}\label{eq:normalized-weak-equation}
 Q(v,\psi)=n\int_Kv^{p-1}\psi\,d\mu.
\end{equation}
\end{lemma}

\begin{proof}
By the definition of $\Lambda$, the functional
$F(w)=Q(w)-\Lambda\|w\|_p^2$ is nonnegative on $\mathcal H$ and
$F(v)=0$. Differentiating $F(v+s\psi)$ at $s=0$, and using
\cite[Theorem~2.6]{LL01} for the  differential of $L^p$ norm, gives
\[
 0=\frac12\left.\frac{d}{ds}\right|_{s=0}F(v+s\psi)
  =Q(v,\psi)-\Lambda\|v\|_p^{2-p}
    \int_Kv^{p-1}\psi\,d\mu.
\]
Since $\Lambda\|v\|_p^{2-p}=\Lambda c^{2-p}=n$, this proves
\eqref{eq:normalized-weak-equation}.
\end{proof}

For every $\varphi\in\mathcal H\cap L^\infty(\mu)$ we can apply 
 \eqref{eq:normalized-weak-equation} to $\psi=v\varphi^2$  and  conclude
\begin{equation}\label{eq:ground-state-identity}
 Q(v\varphi) =Q(v,v\varphi^2)+ \frac4{n-2} \int v^2|D\varphi|^2\,d\mu
 =n\int v^p\varphi^2\,d\mu
 +\frac4{n-2}\int v^2|D\varphi|^2\,d\mu.
\end{equation}

\begin{lemma}\label{lem:all-minimizer-powers}
For every $\beta>0$, the function $v^\beta$ belongs to
$\mathcal H\cap L^\infty(\mu)$. Moreover,
\begin{equation}\label{eq:minimizer-powers}
 D(v^\beta)=\beta v^{\beta-1}Dv\qquad\text{on }\{v>0\},
\end{equation}
and $D(v^\beta)=0$ almost everywhere on $\{v=0\}$.
\end{lemma}

\begin{proof}
First let $0<\beta<\frac12$. For $\varepsilon>0$ we may use
$\psi_\varepsilon=(v+\varepsilon)^{2\beta-1}\in\mathcal H$
in the Euler equation \eqref{eq:normalized-weak-equation}.
Since $2\beta-1<0$, this gives
\[
 \begin{aligned}
 \frac{4(1-2\beta)}{n-2}
 \int (v+\varepsilon)^{2\beta-2}|Dv|^2\,d\mu
 &=n\int (v-v^{p-1})(v+\varepsilon)^{2\beta-1}\,d\mu\\
 &\leq n\int v^{2\beta}\,d\mu<\infty.
 \end{aligned}
\]
Thus $(v+\varepsilon)^\beta$ is bounded in $\mathcal H$.
It converges to $v^\beta$ in $L^2(\mu)$, so weak compactness
implies $v^\beta\in\mathcal H$. For $\beta\geq\frac12$, use
$v^\beta=(v^{1/4})^{4\beta}$ and the usual chain rule.
Finally, the chain rule on $\{v>\delta\}$, followed by
$\delta\downarrow0$, gives \eqref{eq:minimizer-powers};
locality gives the assertion on $\{v=0\}$.
\end{proof}

\section{A weighted variance estimate and the choice of center}
\label{sec10-ce-section}
All integrals in this section are over $K$. We fix a bounded nonnegative
Borel representative of the normalized minimizer $v$  from
Section~\ref{sec:sobolev-minimizer}. Put
\begin{equation}\label{sec10-ce-conformal-notation}
 f=v^{2/(n-2)},\qquad d\eta=f^{n+1}\,d\mu.
\end{equation}
Thus $v^2=f^{n-2}$, $v^p=f^n$   are in $\mathcal H\cap L^{\infty}$ and $W=\int f^n\,d\mu>0$.

\begin{lemma}\label{sec10-ce-center}
The total second moment of $\eta$ satisfies
\begin{equation}\label{sec10-ce-total-variance}
 \iint_{K\times K}d(x,y)^2\,d\eta(x)d\eta(y)\geq2W^2.
\end{equation}
\end{lemma}

\begin{proof}
For a fixed $y\in Y$, put $r=d(y,\cdot)$ and apply
\eqref{eq:radial-variational-sobolev} to $g=f^n=v^p$ to obtain
\[
 W\leq\int r\left(v^p\sqrt{1-|Dr|^2}
                  -\frac2{n-2}v^{p-1}\langle Dr,Dv\rangle\right)\,d\mu.
\]
Using $-\langle Dr,Dv\rangle \leq |Dv||Dr|$ and $a_1a_2+a_3a_4\leq \sqrt {a_1^2+a_3^2} \sqrt{a_2^2+a_4^2}$, we conclude
\[
 v^p \sqrt{1-|Dr|^2}-\frac2{n-2}v^{p-1}\langle Dr,Dv\rangle
 \leq\left(v^{2p}+\frac4{(n-2)^2}v^{2p-2}|Dv|^2\right)^{1/2}.
\]
Integrating this inequality and rewriting the integrand in terms of $f$, we obtain
\[
 W\leq\int r f^{\frac {n+1} 2}
 \left(f^{n-1}+\frac4{(n-2)^2}f|Dv|^2\right)^{1/2}\,d\mu.
\]
An application of Cauchy--Schwarz yields
\begin{equation}  \label{eq:second-integral}
 W^2\leq
 \left(\int r^2f^{n+1}\,d\mu\right)
 \left(\int(f^{n-1}+\frac4{(n-2)^2}f|Dv|^2)\,d\mu\right)\,.
\end{equation}
The second integral can be evaluated using 
\eqref{eq:normalized-weak-equation} for $\psi=vf=v^{p/2}$. In this case \[D\psi= \frac  p 2  \cdot v^{\frac p 2 -1}Dv= \frac n {n-2}fDv\,.\]
Thus, the Euler equation \eqref{eq:normalized-weak-equation}  yields, after division by $n$:
\[
\frac4{(n-2)^2}\int f|Dv|^2\,d\mu
+\int v^2f\,d\mu
=\int v^pf\,d\mu.
\]
Since $v^2=f^{n-2}$ and $v^p=f^n$, the second integral in \eqref{eq:second-integral} equals  $\int f^{n+1}\,d\mu$, hence
\begin{equation} \label{eq:w2}
	W^2\leq  \left(\int r^2f^{n+1}\,d\mu\right)  \left(\int f^{n+1}\,d\mu\right)   =\int r^2\, d\eta \cdot \eta (K) = \eta (K) \int d(x,y)^2\,d\eta (x)\,.
	\end{equation}
This  is true for every $y\in Y$. The CAT(0) variance inequality \cite[Theorem 4.9(iii)]{Sturm}  gives
\[\iint d(x,y)^2\,d\eta(x)\,d\eta(y) \geq 2\eta(K)\inf_{z\in Y}\int d(x,z)^2\,d\eta(x)
\geq 2W^2.
\]
This proves the lemma.
\end{proof}

Let $G\subset K$ consist of points $y$ with $f(y)>0$ which
are Lebesgue points of $f$ and at which $\mu$ has a positive integer
$n$-dimensional density $\theta_y$.   Explicitly,
\begin{equation}\label{sec10-ce-density}
	\lim_{\delta\downarrow0}\frac{\mu(B_\delta(y))}{\omega_n\delta^n}
	=\theta_y,\qquad
	\lim_{\delta\downarrow0}\delta^{-n}
	\int_{B_\delta(y)}|f-f(y)|\,d\mu=0.
\end{equation}
 By Theorem \ref{thm:charts}  and the
Lebesgue differentiation theorem, $f^n\,d\mu$ is concentrated on $G$.
We rewrite the estimate of Lemma~\ref{sec10-ce-center} as 
\[
 \int_K f(y)\left(\int_K d(x,y)^2\,d\eta(x)\right)
               f(y)^n\,d\mu(y)\geq2W^2.
\]
Since  $\int f^n\,d\mu=W$, there is a center $y\in G$ such that
\begin{equation} \label{eq:choice}
 f(y)\int_K d(x,y)^2\,d\eta(x)\geq2W.
\end{equation}
We fix such $y$ until the end of the proof of Theorem~\ref{thm:main}.

\section{Variation and monotonicity of the Sobolev energy}
\label{sec:energy-monotonicity}

We work with $y$ chosen in  Section~\ref{sec10-ce-section} and put
\[
 r=d(y,\cdot),\qquad
 \rho=r\sqrt{f(y)f}.
\]
By Lemma~\ref{lem:all-minimizer-powers},
$\rho=\sqrt{f(y)}\,r\,v^{1/(n-2)}$ belongs to
$\mathcal H\cap L^\infty(\mu)$.
For $a\geq0$, define
\[
 h_a(t)=(1+at^2)^{-(n-2)/2},\qquad
 v_a=vh_a(\rho),\qquad E(a)=Q(v_a).
\]
For $a>0$, the radial functions $x\mapsto h_a(|x|)$ are the
Aubin--Talenti bubbles on $\R^n$, which minimize the Euclidean Sobolev quotient
\cite[p.~353]{Talenti}.

Since $v,\rho\in\mathcal H\cap L^\infty(\mu)$, the map $a\mapsto v_a$ is $C^1$ with values in
$\mathcal H$, with
\begin{equation} \label{eq:pava}
 \partial_a v_a=-\frac{n-2}{2}\rho^2v(1+a\rho^2)^{-n/2},
 \qquad E'(a)=2Q(v_a,\partial_a v_a).
\end{equation}
In particular, $E$ is positive and $C^1$, and $E(0)=nW$.
At $a=0$, the Euler equation gives
\begin{equation}\label{eq:initial-energy-derivative}
 E'(0)=-(n-2)Q(v,v\rho^2)
      =-n(n-2)\int\rho^2v^p\,d\mu
      =-n(n-2)\int r^2f(y)ff^n\,d\mu.
\end{equation}
Hence, the choice of $y$ and \eqref{eq:choice} show
\begin{equation}
E'(0)\leq -2n (n-2) W.
\end{equation} 

To estimate $E(a)$  for large $a$, put
\[
 J(a)=\left(\frac{E(0)}{E(a)}\right)^{2/(n-2)}.
\]

\begin{proposition}\label{sec10-ce-monotonicity}
The function $(J(a)-1)/a$ is nondecreasing for $a>0$.
\end{proposition}

\begin{proof}
Since $E$ is $C^1$ and positive, so is $J$. Differentiating
$(J(a)-1)/a$ and multiplying by $a^2E(a)/J(a)>0$, we see that
it suffices to prove
\begin{equation}\label{sec10-ce-reduced-inequality}
 -\frac{2a}{n-2}E'(a)-E(a)+\frac{E(a)}{J(a)}\geq0.
\end{equation}
We have  $E(a) \geq \Lambda \|v_a\|_p^2$ by the definition of $\Lambda$.    
Since $E(0)=nW$ and $p/2=n/(n-2)$, 
\[	\frac{E(a)}{J(a)}
	=E(0)^{-\frac2{n-2}}E(a)^{\frac n{n-2}}
	=(nW)^{-\frac2{n-2}}E(a)^{\frac p2}
	\geq(nW)^{-\frac2{n-2}}
	\Lambda^{\frac p2}\|v_a\|_p^p
\]
Since $(nW)^{-2/(n-2)} \Lambda ^{p/2}=n$, it remains to show
\begin{equation}\label{sec10-ce-key-inequality}
 -\frac{2a}{n-2}E'(a)-E(a)+n\int v_a^p\,d\mu\geq0.
\end{equation}

For any $g\in \mathcal H\cap L^{\infty} (\mu)$, we apply \eqref{eq:radial-quadratic-sobolev}  to the nonnegative function $g^2\in \mathcal H$ and deduce
\begin{equation}\label{sec10-ce-defect}
\mathcal D(g^2):= \frac n 4 \int r^2g^2\,d\mu  -\int\bigl(ng^2|Dr|^2+2rg\langle Dr,Dg\rangle\bigr)\,d\mu \geq 0.
\end{equation}
On the other hand, for any such $g$,
\[
	Q(gr) =\frac 4 {n-2} \int \left(g^2  |Dr|^2+r^2|Dg|^2+2rg\langle Dr,Dg\rangle\right) \, d\mu + n\int r^2g^2  \, d\mu
\]
Therefore, we can rewrite $Q(rg)-4\mathcal D(g^2)$ as 
\begin{equation} \label{eq:radial-square}
	\begin{aligned}
&Q(rg)-4\mathcal D(g^2)\\		
 =&\frac 4 {n-2}\int  \left((1 +n(n-2)) g^2|Dr|^2  +(2+2(n-2))rg \langle Dr,Dg\rangle + r^2 |Dg|^2\right) \,d\mu \\
=&\frac 4 {n-2}\int \left| rDg + (n-1) gDr\right|^2 \, d\mu
\end{aligned}
\end{equation}

Fix $a>0$.  We are going   to identify the left-hand side of \eqref{sec10-ce-key-inequality}  with $4\mathcal D(g^2)$, where
\[
g=af(y)\,r\,v_a^{p/2}\in\mathcal H\cap L^\infty(\mu).
\]
By the definitions of $\rho$ and $v_a$, and by \eqref{eq:pava}, we have
\begin{equation} \label{eq:rg}
rg=\frac{a\rho^2}{1+a\rho^2}v_a
=-\frac{2a}{n-2}\partial_a v_a,
\qquad
v_a-rg=v(1+a\rho^2)^{-n/2}.
\end{equation}
The first equality together with \eqref{eq:pava} shows
\[-\frac {2 a } {n-2} E'(a)= 2 Q(v_a,rg)
\]
Therefore,   we can rewrite two summands in \eqref{sec10-ce-key-inequality} as
\begin{equation}
	-\frac {2 a } {n-2} E'(a)-E(a)=2Q(v_a,rg) -Q(v_a,v_a) =Q(rg)-Q(v_a-rg)\,.	
\end{equation}
It remains to verify the equality
\begin{equation}
	Q(rg)-Q(v_a-rg) +n\|v_a\|_p^p =4\mathcal D(g^2).	
\end{equation}
Due to \eqref{eq:radial-square}, we only need to show 
\begin{equation}
Q(v_a-rg) -n\|v_a\|_p^p =\frac 4 {n-2}\int \left| rDg + (n-1) gDr\right|^2 \, d\mu
\end{equation}
By \eqref{eq:rg} and  the identity \eqref{eq:ground-state-identity}, applied to the function   $\varphi=(1+a\rho^2)^{-n/2}$ we obtain
\begin{equation} \label{eq:lastintegral}
Q(v_a-rg)=n \int v^p (1+a\rho^2)^{-n} \, d\mu + \frac 4 {n-2} \int v^2 |D\left((1+a\rho^2)^{-\frac n 2}    \right) |^2 \, d\mu
\end{equation}
Since $v^p(1+a\rho^2)^{-n} =v_a^p$, it remains to verify that
\begin{equation} \label{eq:integrands}
	 \int v^2 |D\left((1+a\rho^2)^{-\frac n 2}    \right) |^2 \, d\mu =
	\int
	|rDg+(n-1)gDr|^2\,d\mu.
\end{equation}
Indeed, we are going to prove the stronger almost everywhere pointwise identity:
\begin{equation}  \label{eq:pointwiseid}
-	v D\left((1+a\rho^2)^{-\frac n 2}    \right) =rDg+(n-1)gDr\,.	
\end{equation}
This equality can be verified by expanding all terms.  We are going to provide a slightly less  cumbersome explanation. Both sides are $0$ almost everywhere on $\{v=0\}$.  Thus, we may restrict attention to the subset $\{v>0\}$.   There we have the identity:
\[
	v_a^{-2/(n-2)}
	=v^{-2/(n-2)}+af(y)r^2.
\]
Taking differentials on both sides, we obtain
\[
v_a^{-\frac n {n-2}}Dv_a
=v^{-\frac n {n-2}}Dv-(n-2)af(y)r\,Dr.
\]
Multiplying both sides by $v_a^{n/(n-2)}$, and using the definition of $g$, we see
\begin{equation}  \label{eq:dvafirst}
\begin{aligned}
	Dv_a
	&=\left(\frac{v_a}{v}\right)^{n/(n-2)}Dv
	-(n-2)gDr\\
	&=(1+a\rho^2)^{-n/2}Dv-(n-2)gDr.
\end{aligned}
\end{equation}
On the other hand, due to \eqref{eq:rg},
\[
v_a=rg+v(1+a\rho^2)^{-n/2}.
\]
Differentiating this identity gives
\begin{equation}  \label{eq:dvalast}
Dv_a
=rDg+gDr
+(1+a\rho^2)^{-n/2}Dv
+vD\bigl((1+a\rho^2)^{-n/2}\bigr).
\end{equation}
Comparison of  \eqref{eq:dvafirst} with \eqref{eq:dvalast} proves \eqref{eq:pointwiseid} and, therefore,   the proposition.
\end{proof}

The monotonicity formula gives the following estimate.

\begin{corollary}\label{sec10-ce-energy-corollary}
For every $a\geq0$,
\begin{equation}\label{sec10-ce-energy-bound}
 E(a)\leq nW(1+4a)^{-(n-2)/2}.
\end{equation}
\end{corollary}

\begin{proof}
For $a=0$  we have equality in  \eqref{sec10-ce-energy-bound}.  For $a>0$ we apply Proposition~\ref{sec10-ce-monotonicity} and  \eqref{eq:initial-energy-derivative} and conclude
\[
  \frac {J(a)-1}  {a} \geq J'(0)=-\frac2{n-2}\frac{E'(0)}{E(0)}\geq \frac {2} {n-2} \frac {2n(n-2) W} {nW} =4.
\]
Hence  $J(a) \geq 4a+1$, which implies  \eqref{sec10-ce-energy-bound}. 
\end{proof}

\section{Concentration and completion of the proof}
\label{sec:concentration-completion}

We keep the notation of Section~\ref{sec:energy-monotonicity} and consider the  rescaled functions 
\[
\tilde v_a:=a^{(n-2)/4}v_a.
\]
Then  $Q(\tilde v_a)=a^{(n-2)/2}E(a)$, and 
 Corollary~\ref{sec10-ce-energy-corollary}
shows that $(\tilde v_a)$ satisfies 
\begin{equation} \label{eq:finitecontr1}
	\limsup_{a\to \infty}  Q(\tilde v_a) \leq nW \limsup _{a\to \infty} \left(\frac {a} {1+4a} \right)^{\frac {n-2} 2} =\frac {nW} {2 ^{n-2}}.
\end{equation}


On the other hand, we claim 
\begin{equation}\label{eq:finitecontr2}
 \liminf_{a\to\infty} \|\tilde v_a\|_p ^p
 \geq\frac{\theta_y s_n}{2^n}.
\end{equation}
Fix $R>0$. By the definitions of $\tilde v_a$ and $\rho$,
\begin{equation} \label{eq:intb}
	\int_{B_{R/{\sqrt a}}(y)}\tilde v_a^p\,d\mu
	=a^{n/2}\int_{B_{R/\sqrt a}(y)}
	\left(\frac{f}{1+af(y)r^2f}\right)^n\,d\mu \,.
\end{equation}
The function $f$ is bounded and  the maps $s\mapsto(s/(1+cs))^n$, $c\geq0$,
are uniformly Lipschitz on $[0,\|f\|_\infty]$.  Hence,  for a constant $b>0$, depending only on $\|f\|_{\infty}$ we have
\begin{equation} \label{eq:lpest}
\left|	\left(\frac{f (x)}{1+af(y)r^2(x)f(x)}\right)^n  -	\left(\frac{f (y)}{1+af(y)r^2(x)f(y)}\right)^n \right | \leq b|f(x)-f(y)|\,.
\end{equation}
Since $y$ is a Lebesgue density point of $f$, see   \eqref{sec10-ce-density},
we have 
\[
\lim_{a\to \infty} a^{n/2}
\int_{B_{R/\sqrt a}(y)}|f-f(y)|\,d\mu=0.
\]
Thus, combining \eqref{eq:lpest} and \eqref{eq:intb}, we obtain
\[\liminf  _{a\to \infty}	\int_{B_{R/{\sqrt a}}(y)}\tilde v_a^p\,d\mu =\liminf_{a\to \infty}   a^{\frac n 2} \int_{B_{R/\sqrt a}(y)}
\left(\frac{f (y)}{1+af(y)^2r^2}\right)^n\,d\mu \,.   \]
 
By the  density assertion \eqref{sec10-ce-density}, for every $t>0$, we have
\[
\lim _{a\to \infty} a^{n/2}\mu(B_{t/\sqrt a}(y))
= \theta_y\omega_n t^n.
\]
Integrating the radial function and using this convergence,
we obtain
\[
	\lim_{a\to\infty}
	\int_{B_{R/\sqrt a}(y)}\tilde v_a^p\,d\mu
	=\theta_y n\omega_n f(y)^n
	\int_0^R\frac{t^{n-1}}
	{(1+f(y)^2t^2)^n}\,dt
	=\theta_y n\omega_n
	\int_0^{Rf(y)}\frac{t^{n-1}}
	{(1+t^2)^n}\,dt.
\]

Since $f(y)>0$, letting $R\to\infty$ now gives
\[
	\liminf_{a\to\infty}\|\tilde v_a\|_p^p \geq \theta_y n\omega_n
	\int_0^\infty\frac{t^{n-1}}{(1+t^2)^n}\,dt
	=\frac{\theta_y s_n}{2^n}.
\]
This proves \eqref{eq:finitecontr2}.

By the definition of $\Lambda$ and \eqref{eq:finitecontr2},
\[
	\liminf_{a\to\infty}Q(\tilde v_a)
	\geq
	\Lambda\left(
	\liminf_{a\to\infty}\|\tilde v_a\|_p^p
	\right)^{2/p}
	\geq
	nW^{2/n}
	\left(\frac{\theta_y s_n}{2^n}\right)^{2/p}
	=\frac{nW^{2/n}}{2^{n-2}}
	(\theta_y s_n)^{(n-2)/n}.
\]
Comparing  with \eqref{eq:finitecontr1}, we obtain
\[
W\geq\theta_y s_n\geq s_n,
\]
contrary to $W<s_n$. This completes the inductive
step for $n>2$ and finishes the proof of
Theorem~\ref{thm:main} and
Proposition~\ref{prop:variational}.

\section{Equality and rigidity}
\label{sec:equality-rigidity}

\subsection{Recognition of flat parts}   We begin with two characterizations of flat subsets.

\begin{proof}[Proof of Proposition \ref{rig:variance-flatness}]
	Put $I(y)=\int d(x,y)^2\,d\sigma(x)$. Then
\begin{equation} \label{eq:sturm}
	I(y)\ge I(o)+d(o,y)^2\qquad(y\in Y),
	\end{equation}
	by the variance inequality \cite[Proposition 4.4]{Sturm}.
	Integrating \eqref{eq:sturm}  and using \eqref{rig:variance-flatness-assumption}
	shows  that equality holds in \eqref{eq:sturm} for $\sigma$-almost every $y\in Y$.  By continuity, that equality holds for every $y$  in the support $K$ of $\sigma$.
	
	Recall that a function $f$ is  $\lambda$-convex, respectively $\lambda$-affine on a geodesic $\gamma$ parametrized by arclength, if
	$t\to f\circ \gamma (t)-\frac \lambda 2 t^2$ is convex, respectively, affine  on  the interval of definition of $\gamma$.  Since $Y$ is CAT(0), for every $x\in Y$ the function  $d_x^2$  is $2$-convex on every geodesic in $Y$. 
	
	By integration, the function $I$ is $2$-convex on each geodesic in $Y$.
	Consider the geodesic $\gamma_y:I\to Y$ from the barycenter $o$ to $y$ parametrized  by arclength.  Then $d_o^2$ is $2$-affine on $\gamma _y$.  
	 Equality in \eqref{eq:sturm} implies that $d_x^2$ is  $2$-affine on $\gamma _y$ for almost every $x$.  By continuity, $d_x^2$ is $2$-affine on $\gamma_y$ for every $x\in K$.

	The flat triangle lemma \cite[p. 117]{AKP} shows that for all $x,y \in K$
	the triangle $o,x,y$ is flat.   Hence, the $1$-Lipschitz map  $\log_o:Y\to T_oY$ restricts to an isometric embedding $\ell:K \to \ell (K)=K' \subset T_oY$.
	Denote by $C$ the  closed convex hull of $K'$ in $T_oY$.   
	
	If $C$ is isometric to a convex subset of a Hilbert space, then the inverse isometry $\ell^{-1}:K'\to K$ extends to a  $1$-Lipschitz map $F:C\to Y$ by the Kirszbraun--Lang--Schroeder theorem \cite{LangSchroeder}, \cite[Theorem 10.14]{AKP}.  The composition of $1$-Lipschitz maps
	\[G=P_C\circ \log_o \circ F:C\to C,\] 
	where $P_C$ is the closest-point projection from the CAT(0) space $T_oY$ onto the convex subset $C$, is $1$-Lipschitz and identity on $K'$.  Due to the uniqueness of geodesics, $G$ is the identity on the closed convex hull $C$ of $K'$.  Hence, $F$ is an isometric embedding and $F(C)$ is the closed convex hull of $K$ in $Y$.
	
	It remains to verify that the closed convex hull $C$ of $K'$ in the  cone $T_oY$ is flat.

	Put $|w|=d(0,w)$ for $w\in T_oY$.  For $u\in T_oY$, consider  the function $h_u:T_oY\to \R$
	\[
	h_u(w)=\langle u,w\rangle =\frac12\bigl(|u|^2+|w|^2-d(u,w)^2\bigr).
	\]
	This function $h_u$ is concave on $T_oY$, cf. \cite[Section 3.2]{BuyaloSchroeder}.  By the first variation formula \cite[Theorem 9.36]{AKP}, $-h_{\ell(x)}(w)$ is the derivative of $\frac 1 2 d_x^2$ at the point $o$ in the direction $w$.  Hence, the derivative of $-I$ at $o$
	in the direction $w$ is given as 
	\[H(w):=-D_oI(w)  =2\int_K h_{\ell(x)}(w)\,d\sigma(x).
	\]
	Hence $H$ is concave. Since $o$ is the maximum of $-I$,  $H\leq 0$ on $T_oY$.  For every $y\in K$
	\[
	H(\ell(y))= 2\int _K h_{\ell (x)} ( \ell(y)) \, d\sigma (x)
	= \int _K (d(o,x)^2 +d(o,y)^2 -d(x,y)^2)\, d\sigma (x) =I(o)+d(o,y)^2-I(y)\,,
	\]
	where we have  used that $\ell$ is an isometry.  Equality in \eqref{eq:sturm} implies that $H(\ell(y))=0$.  Thus $H$ vanishes on $K'$ and therefore on the closed convex hull $C$ of $K'$.
	
	Since all functions $h_{\ell(x)}$ are concave, the vanishing of $H$ on $C$ implies that $h_{\ell (x)}$ is affine on $C$ for almost every $x\in K$.  By continuity, $h_u$ is affine on $C$ for all $u$ in $K'$.

	By \cite[Theorem~1.2]{LytchakSchroeder},
	there is an isometric embedding into a direct product
	\[	\iota:C\to \hat C\times \mathbb H,	\]
	where $\mathbb H$ is a Hilbert space, such that every
	Lipschitz affine function on $C$ factors through
	the projection $P_\mathbb H:C\to \mathbb H$, preserving its
	Lipschitz constant.
	
	For $u,v\in K'$, the function $h_u-h_v$ is
	affine on $C$ and is $d(u,v)$-Lipschitz by
	the CAT(0) quadrilateral inequality.
	Moreover,
	\[
	(h_u-h_v)(u)-(h_u-h_v)(v)=d(u,v)^2.
	\]
	The factorization through $P_\mathbb H$ therefore gives
	\[
	d(u,v)^2
	\le d(u,v)\,\|P_\mathbb H(u)-P_\mathbb H(v)\|
	\le d(u,v)^2.
	\]
	Thus $P_\mathbb H$ preserves the distances between
	all points of $K'$.
	
	Since $\iota$ is an isometric embedding into
	the metric product $\hat C\times \mathbb H$, the
	projection of $K'$ to $\hat C$ consists of
	a single point. The corresponding Hilbert
	fiber is closed and convex in $C$.
	Since $C$ is the closed convex hull of $K'$,
	it is contained in this fiber and hence flat.
\end{proof}

The next lemma characterizes boundaries of isometrically embedded Euclidean balls.

\begin{lemma}\label{rig:round-sphere}
	Let $0\ne A\in\I_n(Y)$ be a cycle with compact support $K$,
	put $\mu=\|A\|$, and suppose that $\M(A)=s_n$.
	Let $o$ be the barycenter of $\mu/s_n$ and  assume that 
	\begin{equation}\label{rig:geometric-endpoint-data}
		d(o,x)=1\quad(x\in K),\qquad
		\iint d(x,y)^2\,d\mu(x)d\mu(y)=2s_n^2.
	\end{equation}
	Assume  that, for $\mu$-almost every $x$, the distance function
	$d_x=d(x,\cdot)$ satisfies
	\begin{equation}\label{rig:gradient-endpoint-data}
		|\nabla^A (d_x)|^2=1-\frac{d_x^2}{4}\qquad\mu\text{-almost everywhere}.
	\end{equation}
	Then the closed convex hull of $K$ is a Euclidean unit
	$(n+1)$-ball, and $A$ is its oriented boundary with multiplicity
	one, up to reversal of orientation.
\end{lemma}

\begin{proof}
	Due to \eqref{rig:geometric-endpoint-data} and  Proposition~\ref{rig:variance-flatness}, the closed convex hull of $K$ is flat and we may assume that $Y$ is a Hilbert space $\mathbb H$  and $o$ is  its origin.
	
	Choose a countable dense set $\{x_j\}_{j\geq1}$ of centers for
	which \eqref{rig:gradient-endpoint-data} holds. Since $\mu$ has no
	atoms, we can choose a  point $y\notin\{x_j:j\geq1\}$ at
	which the approximate tangent space  $T_yA$ and   the gradient $\nabla^A (d_{x_j} )_y\in T_yA$ are defined  and   \eqref{rig:gradient-endpoint-data}  holds true  for all $j$.
	
	Since $K$ is contained in the unit sphere $\mathbb S$  of $\mathbb H$, the approximate tangent space $T_yA$ is orthogonal to $y$.  For any other point $x$ in $\mathbb  S$   the $\mathbb H$-gradient $w$ of the distance function $d_x$ is $w=\frac {y-x} {|y-x|}$. The closest point projection  $w'=P_T w$ of $w$  onto $T=T_y\mathbb S$  satisfies
	\[\|w'\|^2 = 1- \frac {\|y-x\|^2} 4 \,.\]
	Thus, the projection $\nabla^A (d_x) _y$ of $w$ onto $T_yA$ satisfies  \eqref{rig:gradient-endpoint-data} if and only if $w'\in T_yA$.  This happens if and only if $x$ is contained in $T_yA \oplus \R\cdot y$.  
	
	Our assumption on $y$ implies that all points $x_j$ are contained in the $(n+1)$-dimensional Euclidean space  $T_yA \oplus \R\cdot y$. Since $x_j$ are dense in $K$ we deduce that all of $K$ is contained in this $(n+1)$-dimensional space, and in fact in its unit sphere $S^n$.
	The constancy theorem for integral currents on the connected
	oriented manifold $S^n$ gives $A=k\llbracket S^n\rrbracket$
	for some nonzero integer $k$. Since $\M(A)=s_n$, we have
	$|k|=1$. In particular, $K=S^n$, and its closed convex hull is
	the Euclidean unit ball.
\end{proof}

\subsection{If-direction}
If $\iota:\overline B^{\,n+1}(0,R)\longrightarrow Y$ is an isometric embedding, then $A=\iota_{\#}\partial[[ B^{n+1}(0,R)]]$ has mass $s_n R^{n}$. Assume that $V(A) < \omega _{n+1} R^{n+1}$. By Proposition \ref{prop:wenger-fillings}, we find  $B\in \I_{n+1} (Y)$ with $\M(B)= V(A)$ and $\partial B=A$.

 The image $Z=\iota (\overline B^{\,n+1}(0,R))$ is convex in $Y$, hence the closest point projection  $P:Y\to Z$ is well-defined and $1$-Lipschitz by the CAT(0) property. Then $B':=P_{\#}(B)$ is a filling of $A$ as well and $\M(B')\leq \M(B)$.  Applying the isometric  inverse $\iota^{-1}:Z\to \R^{n+1}$  we then obtain a filling $\iota^{-1} _{\#} (B')$ of $\partial B^{n+1}(0,R)$   in $\R^{n+1}$  of mass smaller than the volume  of $B^{n+1} (0,R)$.  This contradiction finishes the proof of the if-direction of Theorem \ref{thm:rigidity}.

\subsection{Setting}
\label{rig:analytic-endpoint}
In order to prove the "only if" direction of Theorem \ref{thm:rigidity}, we fix 
 $n\geq1$ and let $0\ne A\in\I_n(Y)$ be a cycle satisfying
\[
 V_Y(A)=c_n\M(A)^{(n+1)/n}.
\]
Put $\mu=\|A\|$ and $K=\spt A$. After rescaling the metric, we may assume that
\begin{equation}\label{rig:normalization}
 m:=\M(A)=s_n,
 \qquad V_Y(A)=\frac{m}{n+1}.
\end{equation}
We are going to find an isometric embedding $\iota$ of the unit Euclidean ball such that $A$ is the push-forward of the boundary of the ball under $\iota$.

\subsection{Dimension one}  \label{subsec:123} Let $n=1$. Write $A=\sum _{i=1}^{\infty}
 A_i$, with $\M(A)= \sum \M(A_i)$,   where every $A_i$ is  the push-forward $(\gamma _i) _{\#}[[\mathbb S^1]]$ for some injective  Lipschitz curve $\gamma _i:\mathbb S^1 \to Y$   \cite[Theorem 5.3]{BDP22}.  Assume that more than one  $A_i$ is nonzero, so we can write $A$ as $A^+ +A^-$ with $\M(A)=\M(A^+)+\M(A^-)$ and
 $A^+,A^-\neq 0$.  Then taking optimal fillings $B^{\pm}$ of $A^{\pm}$ and using Theorem \ref{thm:main} for their masses, we obtain a filling $B=B^++B^-$ of $A$ with 
 \[\M(B)\leq \M(B^+) +\M(B^-) \leq \frac 1 {4\pi} ( \M(A^+) ^2 + \M(A^-)^2)<
 \frac 1 {4\pi} \M(A)^2\,. \]
This contradicts \eqref{rig:normalization} and shows that a non-trivial  decomposition of $A$ is impossible. 
 
 Thus, $A= \gamma_{\#} [[\mathbb S^1]]$ for an injective Lipschitz curve $\gamma$ of length $\ell (\gamma) =\M(A) =2\pi$.  By Reshetnyak's majorization theorem there exists a compact convex domain $\Omega \subset \R^2$ and a $1$-Lipschitz map $\iota:\Omega \to Y$, such that $\iota$ restricted to $\partial \Omega$ is a length-preserving parametrization of $\gamma$. Then $\ell (\partial \Omega ) =2\pi$ and $B=\iota_{\#} [[\Omega]] \in \I_2(Y)$ is  a filling of $A$ of mass at most $\M([[\Omega]])= \mathcal L^2 (\Omega)$.  Thus \eqref{rig:normalization} and the rigidity case of the isoperimetric inequality in $\R^2$  imply that $\Omega$ is isometric to the unit ball, that $\iota$ preserves the mass, and that $B$ is a minimal filling of $A$. 
 
  Since \(\iota\) is \(1\)-Lipschitz and preserves area, the minimality of \(B\) implies, by a cut-and-paste argument, that \(\iota\) is an intrinsic minimal disc in the sense of \cite[Definition~42]{Stadler}. From  \cite[Corollary 48]{Stadler} we conclude that $\iota$ is an isometric embedding.

 \subsection{General preparations}  Assume now $n\geq 2$.  From  the sharp filling inequality for $A$ we derive, for   
   every integral $n$-cycle $B$, 
  \begin{equation}\label{rig:global-variation}
  	c_n\bigl(\M(A)^{\frac{n+1}n}-\M(B)^{\frac {n+1}n }\bigr)\leq V_Y(A-B).
  \end{equation}
 The proof of Lemma \ref{lem:linearized-variation}  provides the following estimate for every $\varepsilon >0$ and 
 \[\delta_{\varepsilon } =c_n m^{\frac {n+1} n} \bigl(1-\bigl (\frac n {n+\varepsilon} \bigr) ^{n+1} \bigr)\,. \]
Every cycle $B\in \I_n (Y)$  with $V_Y(A-B)<\delta_\varepsilon$ satisfies
 \begin{equation}\label{rig:local-variation}
 	\M(A)-\M(B)\leq(n+\varepsilon)V_Y(A-B).
 \end{equation}
 
 The proof of Proposition \ref{prop:compact-support} applies with minor changes of constants  and produces the same conclusions.  In particular, the support $K$  of $A$ is compact. 
 
  Repeating the proof of  Corollary \ref{cor:radial-variation} and using \eqref{rig:local-variation} instead of   \eqref{eq:linearized-variation}  we derive  the following conclusion.  For every $y\in Y$  and the distance function $r$ to the point $y$, for every nonnegative  Lipschitz function $g$ on $K$ and for every $\varepsilon >0$ 
 \[
 \int\bigl(ng+r\langle \nabla ^Ar,\nabla ^Ag\rangle\bigr)\,d\mu
 \leq(n+\varepsilon)\int rg\sqrt{1-|\nabla ^Ar|^2}\,d\mu \,.
 \]
For $\varepsilon\downarrow 0$ we recover the variational inequality \eqref{eq:radial-variational-inequality} for the present cycle $A$.

 The Sobolev calculus developed in
 Section~\ref{sec:sobolev-spaces} applies unchanged.
 Consequently, Corollary~\ref{cor:radial-sobolev}
 remains valid, as does
 Corollary~\ref{cor:intrinsic-sobolev} for $n>2$,
 with exactly the same constants.

\subsection{Dimension two}
 Let $n=2$, so $m=4\pi$. By Proposition~\ref{prop:compact-support}, the lower density of $\mu$ at every  point in $ K$ is at least $1$. Fix a point $y\in K$. Consider the distance function $r$ to $y$ and use the same test functions $g$ as in Section~\ref{sec:dimension-two}. We obtain \eqref{eq:dimension-two-key-estimate} and deduce that the upper density of $\mu$ at $y$ is at most $1$. Thus, the density of $\mu$ at every $y\in K$ is exactly $1$.

Equality in the limit for $\varepsilon\downarrow0$ in \eqref{eq:dimension-two-key-estimate} implies $f-f^2=\frac14$ almost everywhere on $K$, for $f=\frac1r\sqrt{1-|\nabla^A r|^2}$. Hence $f=\frac12$ almost everywhere. This is equivalent to
\begin{equation}\label{rig:dimension-two-gradient}
 |\nabla^A r|^2=1-\frac{r^2}{4}
 \qquad\mu\text{-almost everywhere on }K.
\end{equation}
For every nonnegative Lipschitz function $g$, Corollary \ref{cor:radial-sobolev}   gives the nonnegative expression:
\[
0\leq\mathcal D_y(g)=\int_K\bigl(r^2g-2g-r\langle\nabla^A r,\nabla^A g\rangle\bigr)\,d\mu.
\]
For $0<\varepsilon<1$ consider the function $g_\varepsilon=\max\{r,\varepsilon\}^{-2}$, as in Section~\ref{sec:dimension-two}. Due to \eqref{rig:dimension-two-gradient}, we have  $r\leq 2$, hence $g_{\varepsilon} \geq \frac 1 4$.  Since $\mathcal D_y(g)$ is linear in $g$ and  nonnegative, we have
\begin{equation}\label{eq:1}
 0\leq \frac 1 4 \cdot\mathcal D_y(1)\leq\mathcal D_y(g_\varepsilon)\,.
\end{equation}
 On the other hand, we derive as in Section~\ref{sec:dimension-two}, using  a sequence $\varepsilon\downarrow0$ with $\mu(\{r=\varepsilon\})=0$
\[
 \mathcal D_y(g_\varepsilon)
 =\frac12\mu(\{r>\varepsilon\})
   -\frac{2\mu(B(y,\varepsilon))}{\varepsilon^2}
   +\frac1{\varepsilon^2}\int_{B(y,\varepsilon)}r^2\,d\mu
 \longrightarrow\frac m2-2\pi=0.
\]
Therefore \eqref{eq:1} implies $\mathcal D_y(1)=0$. In the previous notation
\[
 I(y):=\int_Kd(x,y)^2\,d\mu(x)=\int_Kr^2\,d\mu,
\]
this means $I(y)=2m$.
Since $y\in K$ was arbitrary,  this applies to every $y\in K$.

Let $o$ denote the barycenter of $\mu/m$. By the variance inequality \cite{Sturm},
\[
 2m^2=\int_K I(y)\,d\mu(y)\geq 2mI(o),
\]
so $I(o)\leq m$. On the other hand, the radial inequality  \eqref{eq:radial-variational-inequality} with $g=1$ and center $o$, followed by Cauchy--Schwarz, gives
\[
 m\leq\int_Kd_o\sqrt{1-|\nabla^A d_o|^2}\,d\mu
 \leq\int_Kd_o\,d\mu\leq\sqrt{mI(o)}.
\]
Hence $I(o)=m$ and equality holds throughout. In particular, equality in Cauchy--Schwarz implies $d(o,x)=1$ for $\mu$-almost every $x$, and therefore for all $x\in K$ by continuity. Moreover,
\[
 \iint_{K\times K}d(x,y)^2\,d\mu(x)d\mu(y)=2m^2.
\]
Together with \eqref{rig:dimension-two-gradient}, these are the hypotheses of Lemma~\ref{rig:round-sphere}. An application of that Lemma finishes the proof of   Theorem~\ref{thm:rigidity} for $n=2$.

\subsection{Higher dimensions}  
We observe that  the last part of the above proof  was not 
specific to dimension two.  More precisely, after the verification of  \eqref{rig:dimension-two-gradient} and of the equality $I(y)=2m$ for all $y\in K$, the proof works for all dimensions.  Thus, in order to finish the proof of Theorem \ref{thm:rigidity} it suffices to prove the following result, where we use $I(y)=\int d(x,y)^2 \, d\mu (x)$.

\begin{proposition}\label{rig:analytic-identities}
The multiplicity of $A$ is one $\mu$-almost everywhere. $I(y)=2m$, for all $y\in K$. There is a
set $G\subset K$ of full $\mu$-measure such that, for every $y\in G$,
\begin{equation}\label{rig:distance-gradient}
 |Dd_y|^2=1-\frac{d_y^2}{4}
 \qquad\mu\text{-almost everywhere on }K.
\end{equation}
\end{proposition}

In order to prove the proposition, we may and will assume from now on $n\geq 3$.  We are going to use notation for $p, k_n, Q, \Lambda$ introduced in Sections \ref{sec:intrinsic-sobolev} and  \ref{sec:sobolev-minimizer}. We have $0<\Lambda\leq k_n$.  If $\Lambda<k_n$, the
minimizer-existence argument of Proposition~\ref{prop:normalized-minimizer}
still applies: the strict inequality $\Lambda<k_n$ replaces the
assumption $m<s_n$ used there.  The rescaled minimizer then satisfies
$W=(\Lambda/n)^{n/2}<s_n$, and the arguments of
Sections~\ref{sec10-ce-section}--\ref{sec:concentration-completion}
give the same contradiction.  Thus $\Lambda=k_n$.
Since $m=s_n$, the constant function $v=1$ minimizes the quotient
\eqref{eq:spherical-quotient}, and $W=m=s_n$.

The proof of  Lemma \ref{sec10-ce-center}  and of \eqref{eq:choice} works without changes.  If we do not have equality in \eqref{eq:choice} for almost all $y$, then we find some $y$ satisfying the conditions listed at the end of Section \ref{sec10-ce-section} and satisfying a strict inequality in \eqref{eq:choice}.  In this case, the proof in the two subsequent Sections lead to a contradiction in the same way.  This already verifies $I(y)=2m$ for almost all $y\in K$.  By continuity, $I(y)=2m$ is true for all $y$ in the support $K$.

  Now all arguments  in Sections \ref{sec:energy-monotonicity}, \ref{sec:concentration-completion} work for almost every center $y\in K$.
  If we do not reach a contradiction  after  \eqref{eq:finitecontr2},  the density $\theta_y$ has to be $1$ at almost every $y\in K$.

 Finally, if the non-decreasing function  $(J(a)-1)/a$  is not constant, then we obtain the strengthening  $E(a)\leq nW (1+ ea)^{-(n-2)/2}$  in \eqref{sec10-ce-energy-bound}, for some $e>4$ and all sufficiently large $a$. Then Section \ref{sec:concentration-completion} provides a contradiction.  Therefore, we conclude $J(a)=1+4a$  for all $a\geq 0$ (with respect to almost all centers $y\in K$).
 
 From the proof of Proposition \ref{sec10-ce-monotonicity} we deduce that the nonnegative quantities $\mathcal D(g^2)$   appearing there  vanish for all $a$.   In the present situation, with $v=1$, we have  $g^2=g_1^2=r^2
 (1+r^2)^{-n}$.
We recall the square-completion
argument from the proof of Corollary~\ref{cor:radial-sobolev}:
\[
 \mathcal D(g^2)
 \geq n\int g^2\left(\frac{r^2}{4}+1-|Dr|^2
                    -r\sqrt{1-|Dr|^2}\right)\,d\mu
 =n\int g^2\left(\sqrt{1-|Dr|^2}-\frac r2\right)^2\,d\mu.
\]
Since $0=\mathcal D(g^2)$,
$g>0$ on $K\setminus\{y\}$ and $\mu(\{y\})=0$,
we obtain \eqref{rig:distance-gradient}.

This finishes the proof of Proposition \ref{rig:analytic-identities} and therefore of Theorem \ref{thm:rigidity}.

\end{document}